\documentclass[final]{siamltex}

\usepackage{amsmath}
\usepackage{amssymb}
\usepackage{amstext}

\usepackage{graphicx}
\usepackage{hhline}
\usepackage{psfrag}
\usepackage{float}
\usepackage[dvips]{epsfig}
\usepackage{subfig}
\usepackage{tikz}

\begin{document}

\title{A meshfree method for the BGK model for rarefied gas dynamics}

\author{ S. Tiwari  \footnotemark[1] , 
A. Klar \footnotemark[1] \footnotemark[2]  \and G. Russo \footnotemark[3] }
\footnotetext[1]{Technische Universit\"at Kaiserslautern, Department of Mathematics, Erwin-Schr\"odinger-Stra{\ss}e, 67663 Kaiserslautern, Germany 
  (\{klar, tiwari\}@mathematik.uni-kl.de)}
\footnotetext[2]{Fraunhofer ITWM, Fraunhoferplatz 1, 67663 Kaiserslautern, Germany} 
\footnotetext[3]{Department of Mathematics and Computer Science, University of Catania, Italy (russo@dmi.unict.it)}


\maketitle

\begin{abstract}
 In this paper we have applied a Semi-Lagrangian schemes with meshfree interpolation, based on a Moving Least Squares (MLS) method, to solve the BGK model for rarefied gas dynamics. Sod's shock tube problems are presented for a large range of mean free paths in one dimensional physical space and three dimensional velocity space. In order to validate the solutions obtained from the meshfree method, we have used  the piecewise linear spline interpolation. Furthermore, we have compared the solutions of the BGK model with the solutions obtained from Direct Simulation Monte Carlo (DSMC) method. In the case of a very small mean free path the numerical solutions are compared with the exact solutions of the compressible Euler equations. Overall we found that the meshfree interpolation gives better approximation than the piecewise linear spline interpolation. 
\end{abstract}




\section{Introduction}


The Boltzmann equation is an evolution equation of a probability distribution function consisting  of transport and collision terms \cite{CIP}. Due to the high dimensional integral in the collision term, deterministic 
numerical approaches are  complicated and time consuming. Therefore,  stochastic numerical methods like DSMC, see \cite{Bird, Babovsky, NS},
have been used extensively for complex applications. DSMC methods are suitable for high Mach number and stationary flows. However, for low Mach number flows, the statistical noise inherent in these methods dominates the flow quantities. Here, our main interest is to develop a numerical method for such low Mach number, time dependent flows in arbitrary geometries. Since last 20 years micro-nano scale rarefied gas flows have attracted many researchers due to the fabricated techniques in  
Micro-Electro-Mechanical-Systems (MEMS) \cite{Gad}, devices, for examples, micro pump, micro turbines, micro pipes \cite{KBA}.  We consider a simplified model for rarefied gas flows, where  deterministic methods can be applied more easily. We choose a simplified model suggested by Bhatnager, Gross and Krook \cite{GBK}, the so called    BGK model for the Boltzmann equation, where the collision term is replaced by a relaxation of the distribution function towards a local thermal equilibrium. For  deterministic schemes for this model we refer to \cite{Mieuss} and other references therein. In the present paper, we apply the Semi-Lagrangian scheme suggested by Russo and Filbet, see \cite{RF} for details. In contrast to  
\cite{RF}, where one dimensional physical and one dimensional velocity spaces are considered,  we  consider here  a three dimensional velocity space. Moreover, 
the reconstruction procedure is different compared to the one applied in \cite{RF}. Here, we use a meshfree method for the reconstruction. Meshfree methods are suitable for changing computational domains in time or flows in 
complicated geometries, see \cite{TK07, TKH09}. We note that a  meshfree method based on  Least-squares was applied to solve the compressible Euler equations, see \cite{Deshpandeetal} 
and other references there. 

The paper is organized as follows. In section \ref{sec:model} the BGK model for the Boltzmann equation is presented. In section \ref{sec:num_scheme} the semi Lagrangian scheme for the model and the 
boundary conditions are described. In section \ref{sec:interpolation} we present piecewise linear spline interpolation and the moving least squares (MLS) approximation for the reconstruction of the function. 
In section \ref{sec:numerics} Sod's shock tube problem \cite{SOD} is solved for several range of mean free paths. For larger Knudsen numbers the numerical solutions for the BGK model obtained from the 
piecewise linear spline and MLS interpolations are compared with the solutions obtained from the DSMC simulations of the Boltzmann equation. For a very small Knudsen number, we have compared the 
numerical solutions  of the BGK model with the exact solutions of the compressible Euler equations. We found that the solutions obtained from the MLS approximation are closer to the DSMC results or the  exact solutions than the solutions 
obtained from the piecewise linear spline interpolation. Moreover, we have compared  numerical approaches based on continuous and discrete Maxwellians as suggested in \cite{Mieuss}. 
We found that the use of discrete Maxwellian allows us to reduce the number of velocity grids, which is very important in higher dimensional cases from the memory as well as computational point of view. 
Finally, in section \ref{sec:conclusion} some conclusions and future works are presented.

\section{The BGK model for rarefied gas dynamics}
 \label{sec:model}
    
The BGK model is the simplified model of the Boltzmann equation for a rarefied gas dynamics, where the collision term is modeled by a relaxation of the distribution function $f(t,x,v)$ to the Maxwelian equilibrium distribution. The is the evolution equation of the distribution function $f(t,x,v)$ and is given by the following initial boundary value problem
\begin{equation}
\frac{\partial f}{\partial t} + v_x \frac{\partial f}{\partial x} = \frac{1}{\tau}(M - f)
\label{bgk_eqn}
\end{equation}
with $f(0,x, v) = f_0(x, v), \; t \ge 0, x \in [a, b] \subset \mathbb{R}, \;{ v} \in \mathbb{R}^3$  and some boundary conditions assigned at $a$ and $b$, as will be described in the next section. We denote by $v = (v_x, v_y, v_z)$ the $3d$ velocity vector. 
 Here $\tau$ is the relaxation time and $M$ is the local Maxwellian given by 
\begin{equation}
M = \frac{\rho}{(2\pi R T)^{3/2}} \exp \left (-\frac{| v -  U|^2}{2RT}\right), 
\label{maxwellian}
\end{equation}
where the parameters $\rho, U, T$ are macroscopic quantities, namely, density, mean velocity and temperature, respectively. Here, $R$ is the gas constant. The macroscopic quantities $\rho, U ,T $ are computed from $f(t, x,  v)$ as its moments. 
In this case we have denoted $U = (U_x, U_y, U_z)$. 
Let $ { \phi } (v)=\left (1,  v ,\frac{| v|^2}{2} \right )$ be the collision invariants. The moments are defined by 
\begin{equation}
(\rho, \rho U, E) = \int_{\mathbb{R}^{3}} \phi(v) f(t,  x,  v) d v.
\label{moments}
\end{equation}
Here, $E$ is the total energy density and it is related to the temperature through the internal energy 
\begin{equation}
e(t, { x}) = \frac{3}{2}R T, \quad \quad \rho e = E - \frac{1}{2}\rho |U|^2. 
\label{internal_energy}
\end{equation}
The relaxation time $\tau$ and the mean free path $\lambda$ are related according to  \cite{CC}
\begin{equation}
\tau = \frac{4 \lambda}{\pi \bar C},
\label{tau}
\end{equation}
where $\bar C = \sqrt{\frac{8RT}{\pi}}$ and the mean free path is given by 
\[
\lambda = \frac{k_B}{\sqrt{2\pi\rho R d^2}},
\]
where $k_B$ is the boltzmann constant and $d$ is the diameter of the gas.  

 \section{Semi-Lagrangian scheme for the BGK model}
\label{sec:num_scheme}

To solve the BGK model, we have used the Semi Lagrangian method suggested by Russo \& Filbet, see  \cite{RF} for details. This method is Semi Lagrangian for the advection and  implicit in the treatment of collision. In this paper we give a short description of the method. 
 We consider constant time step $\Delta t$, uniform meshes in velocity space  with mesh size $\Delta v$ and in physical space not necessarily uniform meshes with average spacing $\Delta x$. Let $t_{\rm final}$ be 
 the final time step of computation. The time steps are  given by  $t_n = n\Delta t, n = 0, 1, \ldots $. The space discretization is obtained by generating grid 
 points (regular or irregular) ${ x}_i  \in [a, b], i = 1 , \ldots , N_x+1$, where $N_x + 1$ is the total number of grid points in physical space. We note that the $N_x +1 $ grid points include interior as well boundary points $x_1 = a$ and $x_{N_x + 1} = b$. 
Consider the $N_v$ velocity grid points in each directions, where the uniform velocity grid size is given by  $\Delta v = \frac{2  v _{\rm max}}{N_v}$. 
The x-component of velocity grids are defined by  $v _j = -v_{\rm max} + (j-1)\Delta v, j = 1,\ldots , N_v+1$.  Similarly, the $y$- and $z$- components are defined by $v_k$ and $v_l$ for $k,l =1, \ldots, N_v + 1$.  Assuming that $f$ is negligibly small for $|v |>  v_{\rm max}$. 

Let $f_{jkl} = f_{jkl} (t,  x ) = f(t,  x,  v_j, v_k, v_l)$.  The evolution equation of $f_{jkl}(t,  x )$ along the characteristics between time steps $n$ and $n+1$ is calculated from  the Lagrangian form of the discrete BGK model 
\begin{eqnarray}
\frac{df_{jkl}}{dt} &=& \frac{1}{\tau}(M_{jkl}-f_{jkl}) 
\label{bgk_lagrangian}
\\
 \frac{d x}{dt} &= &v_j
\label{bgk_char}
\end{eqnarray}
with initial conditions
\begin{equation}
 x(t_n) = \tilde  x,  \; f_{jkl}(t_n) = f_{jkl}^n(\tilde  x ) =  \tilde{f}_{jkl}^n,  \; t\in[t_n, t_{n+1}]
\end{equation}
together with boundary conditions for $f_{jkl}$ at boundary points. 

Here $M_{jkl}$ is still the local Maxwellian having the same moments of $f_{jkl}$ and is re-expressed by 
\begin{equation}
M_{jkl} = \frac{\rho}{(2\pi R T)^{3/2}} \exp\left(-\frac{(v_j -  U_x)^2 + (v_k-U_y)^2 + (v_l-U_z)^2 }{2RT}\right).
\label{discr_maxw}
\end{equation}

\begin{figure}
 
\begin{center}
	\begin{tikzpicture}
	\begin{scope}[xshift=2cm]    
	\draw[->] (-1,0) node[left]{$t_{n}$} --(7,0) node[right]{$x$};
	\draw[->] (0,-0.2)--(0,3) node[above]{$t$};
	\draw (-1,2) node[left]{$t_{n+1}$}  -- (6.5,2);
	\draw [dashed] (2.5,0)    -- (6,2) node[midway,left]{$v_j > 0$};
          \draw (2,0)  -- (2,2);	 	
         \draw (4,0)   -- (4,2) ;
	\draw (6,0) node [below ]{$x_i^n$} -- (6,2) node[above]{$f_{ijkl}^{n+1}$} node[below]{$x_i^{n+1}$};	
	\draw[->,thick] (1,0.5) node[left]{$\tilde{f}_{ijkl}^n$} -- (2.5,0); 	
         \draw [->,thick] (2,-0.5)node[below]{$\tilde{x}_{ijkl} = x_i^n - v_j \Delta t$} -- (2.5,0); 
         \draw [->,thick] (6.5,0.5) node[right]{$f_{ijkl}^n$} -- (6,0);
 	\end{scope}
	\end{tikzpicture}
	\caption{Computational grid points in space and time} 
	 \label{sketch_char}
\end{center}
\end{figure}
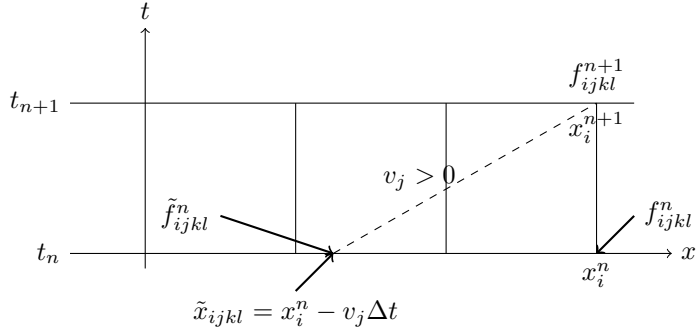

We solve Eq. (\ref{bgk_lagrangian}) by the implicit Euler scheme
\begin{eqnarray}
 f_{ijkl}^{n+1} &= & \tilde{f}_{ijkl}^n + \frac{\Delta t}{\tau_i^n}(M_{ijkl}^{n+1} - f_{ijkl}^{n+1}),
 \label{implicit_euler}
\end{eqnarray}
where $\tau_i^n$ is the relaxation time in grid point $i$ at time level $n$ and the characteristic equation (\ref{bgk_char}) is solved by 
\begin{eqnarray}
 x_i^{n+1} &=& \tilde{  x }_{ijkl} +  v_j \Delta t, \;\; \mbox{for} \;\;i = 1, \ldots , N_x + 1, \; j,k,l = 1,  \ldots , N_v + 1,
\end{eqnarray}
where the initial position $\tilde{x}_{ijkl}$ is given by $\tilde{x}_{ijkl} = x_i-v_j\Delta t$. In Figure \ref{sketch_char} we have given the geometrical interpretation. 
At the time level $t^n$ all values $f^n_{ijkl}, i = 1, \ldots, N_x + 1$ are known. At the time level $t^{n+1}$ the corresponding values are $f_{ijkl}^{n+1}, i = 1, \ldots, N_x + 1$. 

The method consists of three steps: \\
({\bf i})
First, we determine $\tilde{ x}_{ijkl}$ form the backward characteristics  $\tilde{ x }_{ijkl} =  x_i ^{n+1} -  v_j\Delta t$, see Figure \ref{sketch_char}.   Then, we reconstruct (or interpolate)  the function $\tilde{f}^n_{ijkl}$ at $\tilde{x }_{ijkl}$ from the values of its neighboring grid points. 
 One can use any reconstruction, for example, spline interpolations, least squares interpolations.  In this paper we use the piecewise linear interpolation (linear spline) and the linear moving least squares interpolation. Higher order interpolations are also possible, for examples, piecewise cubic spline polynomial \cite{RF}, higher order MLS. 
 Higher order reconstructions give oscillations if the solutions develop shocks, therefore, WENO reconstruction \cite{CV} are necessary to damp the oscillations.

({\bf ii}) In the second step we obtain $M_{ijkl}^{n+1}$.  Since  $M_i^{n+1}$ and $f_i^{n+1}$ have same conservative moments, we multiply the above discrete equation (\ref{implicit_euler}) by the 
collisional invariants $ \phi (  v )$ and sum over the velocity getting
\begin{eqnarray}
\rho_i^{n+1} &=& \sum_{j=1}^{N_v+1}  \sum_{k=1}^{N_v+1}\sum_{l=1}^{N_v+1} \tilde{f}_{ijkl}^n \Delta v^3 
\label{discrete_moments1}
\\
(\rho U_x )^{n+1} &=& \sum_{j=1}^{N_v+1}  \sum_{k=1}^{N_v+1} \sum_{l=1}^{N_v+1} v_{j} \tilde{f}_{ijkl}^n \Delta v^3 
\\
(\rho U_y )^{n+1} &=& \sum_{j=1}^{N_v+1}  \sum_{k=1}^{N_v+1} \sum_{l=1}^{N_v+1} v_{k} \tilde{f}_{ijkl}^n \Delta v^3 
\\
(\rho U_z )^{n+1} &=& \sum_{j=1}^{N_v+1}  \sum_{k=1}^{N_v+1} \sum_{l=1}^{N_v+1} v_{l} \tilde{f}_{ijkl}^n \Delta v^3 
\\
 E_i^{n+1} &=& \frac{1}{2}\sum_{j=1}^{N_v+1} \sum_{k=1}^{N_v+1}\sum_{l=1}^{N_v+1}  (v_j^2 + v_k^2 + v_l^2) \tilde{f}_{ijkl}^n \Delta v^3.
\label{discrete_moments5}
\end{eqnarray}
Now from (\ref{discrete_moments1} - \ref{discrete_moments5}) together with (\ref{internal_energy}) we obtain all five parameters of the Maxwellian and can define $M^{n+1}_{ijkl}$ from (\ref{discr_maxw}). 
\\
({\bf iii}) Finally, we update the density function solving Eq. (\ref{implicit_euler})
\begin{equation}
f_{ijkl}^{n+1} = \frac{\tau^n_i \tilde{f}_{ijkl}^n + \Delta t M_{ijkl}^{n+1}}{\tau^n_i + \Delta t} \;\; \mbox{for} \;\; i = 1, \ldots , N_x + 1, \quad j,k, l = 1, \ldots , N_v + 1. 
\label{SL_step3}
\end{equation}

\subsection{Boundary conditions} We use the  diffuse reflection boundary conditions. This means, when gas molecules hit boundaries, we forget their history. We reflect them according to the 
half Maxwellian with the wall density $\rho_w$,  wall temperature $T_w$ and wall velocity $U_w$. This is the Maxwell boundary condition, see \cite{CIP} for details.  Let $\nu$ be the unit normal on the wall pointing towards the computational domain. In the case of Spline interpolation we need ghost points next to the boundary points $i$ in order to compute $\tilde f^n_{ijkl}$. In order to apply the boundary condition, we first compute $f^{n+1}_{ijkl}$ in all interior points $i$, then we extrapolate the new distribution function $f^{n+1}_{\Gamma}$ for 
$(v-U_w)\cdot \nu > 0$ on the boundary points with the help of MLS interpolation. For $(v-U_w)\cdot \nu < 0$  the diffuse reflection boundary conditions are obtained  according to 
 \begin{equation}
M_{\Gamma}^{n+1} = \frac{\rho_w} {(2\pi R T_w)^{3/2}} \exp\left(-\frac{|v-U_w|^2}{2RT_w}\right)  ,
\end{equation}
where
\[ \rho_w = -\frac{\int_{( { v} - { U}_w)\cdot { \nu} < 0} ( { v} - { U}_w)\cdot { \nu} ~ {f}_{\Gamma}^{n+1} d { v}} { \int_{( { v} - {{U}_w)\cdot { \nu} > 0} } ( { v} - { U}_w)\cdot { \nu} ~ \frac{1} {(2\pi R T_w)^{3/2}} e^{-\frac{|v-U_w|^2}{2RT_w}}  d { v}  } .
\]
 
\section{Interpolation methods}
\label{sec:interpolation}

As we have already mentioned, our main aim is to develop a method to simulate the interactions between rigid body motion and rarefied gas. Due to the movement of a rigid body the computational domain for a gas changes. Moreover, 
the intersection of the surface of a rigid body and cells of rarefied gas makes the numerical scheme more complicated.  In the vicinity of a moving rigid body the regular grid structure does not exists any more. In this section we present two interpolation methods, which are suitable for irregular grids. 
In this paper we consider linear interpolation. Higher order interpolations require some stable reconstructions, like WENO, which will be focused in future works.  

\subsection{Piecewise linear interpolation (or linear Spline $S_1$)}
This is simple to implement and faster than the MLS method. However, one has to add the ghost points to apply boundary conditions, which could be complicated for complex boundaries. 
Let $I_k = [x_k, x_{k+1}]  \subset [a, b]$ be an arbitrary interval and $\tilde x \in [x_k, x_{k+1}] $ be an arbitrary point. The corresponding function values are $f_{k} = f(x_k)$ and $f_{k+1} = f(x_{k+1})$. The linear interpolantion at $\tilde x$ is given 
by 
\begin{equation}
f(\tilde x) = f_k + \frac{f_{k+1}-f_k}{x_{k+1}-x_k}(\tilde x-x_k) = \frac{f_k (x_{k+1}-\tilde x) + f_{k+1}(\tilde x-x_k)}{x_{k+1}-x_k}.
\end{equation}
We note that the size of intervals $I_k$ need not to be equal.

\subsection{Moving least squares (MLS) interpolation}

In contrast to $S_1$ interpolation, this is a fully meshfree method. In the $S_1$ interpolation, only the next left and right grid points are used to interpolate. However, in the MLS approximation, we use the 
nearest neighbor points inside a radius, which is about 2.5 times the average grid space.  Therefore, the computational costs increases slightly in the case of MLS approximation compared to $S_1$ 
interpolation. But the MLS gives better approximation than the $S_1$ interpolation. In this case also the distribution of grids need not to be uniform.
Another advantage of this method is that it is not required to add the ghost points to apply the boundary conditions and is easy to handle complex geometries. 

Let $\tilde x\in [a, b]$ be an arbitrary point. 
 We consider the problem to approximate the function $f = f(\tilde x)$ at $x$ from the values of its neighboring points. We associate a weight function such that the near particles have higher and the far particles have lower influence. In order to limit the number of points the neighboring points are taken those points inside the circle of radius $h$ with 
center $\tilde x$. We choose the radius $h$, for example, some factor of $\Delta x$,  such that we have at least minimum number of neighbors for the least squares approximation. 
Let $P(\tilde x)=\{x_j,j=1, \ldots, m \}$ be the set of $m$ neighbor points of $\tilde x$ inside the radius $h$. 
 We note that this neighboring list is similar to the central stencils in the sense of the finite difference method. Therefore, if the relaxation time $\tau$ is very small and the solution of the Boltzmann equation develop shocks, we need to sort out the neighbor list according to the 
sign of the velocity $v_x$. 
The weight function can be quite arbitrary, but in our computations, we consider a Gaussian weight function 
\begin{eqnarray} 
w_j = w( {x}_j - \tilde x; h) =
\left\{ 
\begin{array}{l}  
\exp \left (- \alpha \frac{(x_j - \tilde x)^2}{h^2} \right), 
\quad \mbox{if    }  \frac{ |x_j - \tilde x|} {h} \le 1 
\\ \nonumber
 0,  \qquad \qquad \mbox{else},
\end{array}
\right.
\label{weight}
\end{eqnarray}
with $ \alpha $ a  user defined positive constant. In our computation, we have considered $\alpha = 6$. Let us sort out the neighboring points from $1$ to $m$ with respect to distance. 
This means, the neighbor index $1$ is the nearest neighbor of $\tilde x$. 

In order to approximate the function we consider the $m$ Taylor's expansions of $f(x_j)$ around $\tilde x$ 
\begin{equation}
f(x_j) = f(\tilde x) + (x_j-\tilde x) \frac{\partial f}{\partial x} + e_j, 
\label{taylor}
\end{equation}
for $j=1, \ldots, m$, where $e_j$ is the error in the Taylor's expansion.  We first assume that $f$ approximates the nearest point $f_1$. In other words,  $e_1 = 0$. 
The unknowns $f, \frac{\partial  f}{\partial x}$ are computed by minimizing the error $e_j$ for 
$j=2, \ldots, m$ and setting the constraint $e_1$ = 0.  To solve this constraint least squares problem, we subtract the first equation with $e_1 = 0$ to all the other equations and the system of equations can be rewritten in the form 
\begin{eqnarray}
f_2 - f_1 &=& (x_2 - x_1) \frac{\partial f}{\partial x} + e_2 \nonumber \\
\vdots & = & \vdots \\
f_m - f_1 &=& (x_m - x_1) \frac{\partial f}{\partial x} + e_m  \nonumber
\end{eqnarray}
The system of equations can be written in the vector form as 
\begin{equation}
{ e } = { b} - M \frac{\partial  f}{\partial x},
\end{equation}
where ${e} = [e_2, \ldots, e_m]^T$,  $b = [f_2 - f_1, \ldots, f_m - f_1]^T $ 
and $ M  = [ x_2 - x_1, \ldots, x_m-x_1]^T $ .
 For $m > 2$, this system of equations is over-determined for one unknown $ \frac{\partial  f}{\partial x}$ . The unknown
$\frac{\partial  f}{\partial x}$ is obtained from the weighted least squares method by minimizing the quadratic form 
\begin{equation}
J = \sum_{j=2}^m w_j e_j^2 = (M \frac{\partial  f}{\partial x}  - { b})^T W (M \frac{\partial  f}{\partial x} - { b}),
\end{equation}
where $W=w_j\delta_{j k}, k = 2,\ldots, m $ is the diagonal matrix. The minimization of $J$ formally yields 
\begin{equation}
\frac{\partial  f}{\partial x} = (M^TWM)^{-1}(M^TW){ b} = \frac{\sum_{j=2}^m w_j (x_j-x_1) (f_j - f_1)}{\sum_{j=2}^m w_j (x_j-x_1)^2}. 
\end{equation}
Now from the equation (\ref{taylor}) with $e_1=0$ for the closest point $x_1$ we can compute the value of $f(\tilde x)$ as 
\begin{equation}
f(\tilde x) = f(x_1) + (\tilde x - x_1) \frac{\partial f}{\partial x} 
\end{equation}
since $\frac{\partial f}{\partial x}$ is now known. 
We note that the higher order approximations are straightforward. Moreover, the approximation in two and three dimensional physical space is also straight forward. 
For example, for higher order, say order $p$, one needs to use the Taylor expansion up to $p$. Then one obtains an overdetermined system of $m-1$ equations in $ p < m$ unknowns which are the derivatives of order $k, k = 1, \ldots, p$. Such a system can be solved in the least squares sense with the technique shown before. In several dimensions, say dimension $d$,  one has to use the Taylor expansion in dimension $d$ up to order $p$, obtaining again an overdetermined set of equations that can be solved in the least squares sense. 
We refer to our earlier papers \cite{TK07, TKH09} for higher orders and several dimensions.  
 
\section{Numerical results}

\label{sec:numerics}
We have considered the Sod's shock tube problem \cite{SOD} as Benchmark to validate our numerical methods. We consider the computational domain $[a, b] = [0, 1]$.  
The initial conditions  are 
\[
\rho^0 = \rho_l, \;\;U^0 = 0, \;\;e^0 = 2.5 \; \; 
\mbox{for} \;\; 0\le x < 0.5
\]
\[
\rho^0 = \rho_r,\;\; U^0= 0,\;\; e^0 = 2.0 \;\; 
\mbox{for}\;\; 0.5\le x \le 1, 
\]
where $\rho_l$ and $\rho_r$ denote the density on the left and right half intervals, respectively. 
We consider the following  boundary conditions 
\[
U(t) = 0, \; e(t) = 2.5 \;\; \mbox{at}\;\; x = 0 \;\;\mbox{and}\;\;
U(t) = 0, \; e(t) = 2 \;\; \mbox{at}\;\; x = 1.
\]
In the MLS we set $\alpha = 6.0$ in (\ref{weight}) and $h$ equal to 
$2.5$ times the initial spacing of the grids. 
The initial spacing of the grids is given by $ \Delta x = 1/N_x$. 
We consider the Argon gas with diameter $d=0.368\times 10^{-9} m$, Boltzmann constant $k_B=1.3806\times 10^{-23} J~K^{-1}$ and the 
 gas constant  $R=208 J {\rm Kg}^{-1}~K^{-1}$. 
The corresponding initial temperature are $0.008012~K$ on the left half of the domain and $0.00641~K$ on the right half of domain. 
The limit of the velocities in all direction is set by $ v_{\rm max} = 10 m/s$.  We note that the thermal velocity corresponding to such a temperature 
is much smaller than the chosen $v_{\rm max}$. 
So, initially the gas is distributed according to the Maxwellian with these initial parameters. 
The final time is $t_{\rm final} = 0.17$ seconds. 

\subsection{Test 1: Comparison of solutions with different CFL numbers}

\begin{figure}
 \begin{center}
\includegraphics[keepaspectratio=true, width=.6\textwidth]{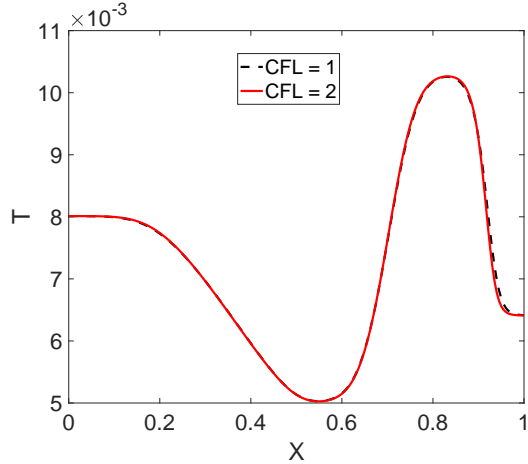}
\caption{Temperature obtained from CFL = 1 and 2 for initial mean free path  ${\lambda}_l=0.001$ m for $x < 0.5$, $\lambda_r = 0.008$ m for $x > 0.5$ with $N_x = 200$ at $t_{\rm final} = 0.17$. }
\label{CFL_plot}
\end{center}
\end{figure}

One of the main advantages of this scheme is that the Courant number CFL $= v_{\rm max}\Delta t/\Delta x$  can be bigger than 1. In the first example we consider $N_x = 200$, $\rho_l = 10^{-4}\, {\rm Kg}\, m^{-1}, \rho_r = 0.125\times \rho_l$, which corresponds to the initial mean free path ${\lambda}_l = 0.001$ m for $x < 0.5$ and ${\lambda}_r = 0.008$ m for $x > 0.5$ and the corresponding initial relaxation time can be computed from the Eq. (\ref{tau}). 
The initial condition is the Maxwellian distribution (\ref{maxwellian}) with parameters given by the initial density, mean velocity and the temperature. 
We have plotted the temperature for CFL = 1 and 2, where the solutions are the same for both cases, see Figure \ref{CFL_plot}. 

For the optimal choice of the CFL number we refer \cite{GGS}. In the following test cases, we have considered CFL = 1. 

\subsection{Test 2: Comparison of solutions with constant and variable $\tau$}

\begin{figure}
 \begin{center}
\includegraphics[keepaspectratio=true, width=.3\textwidth]{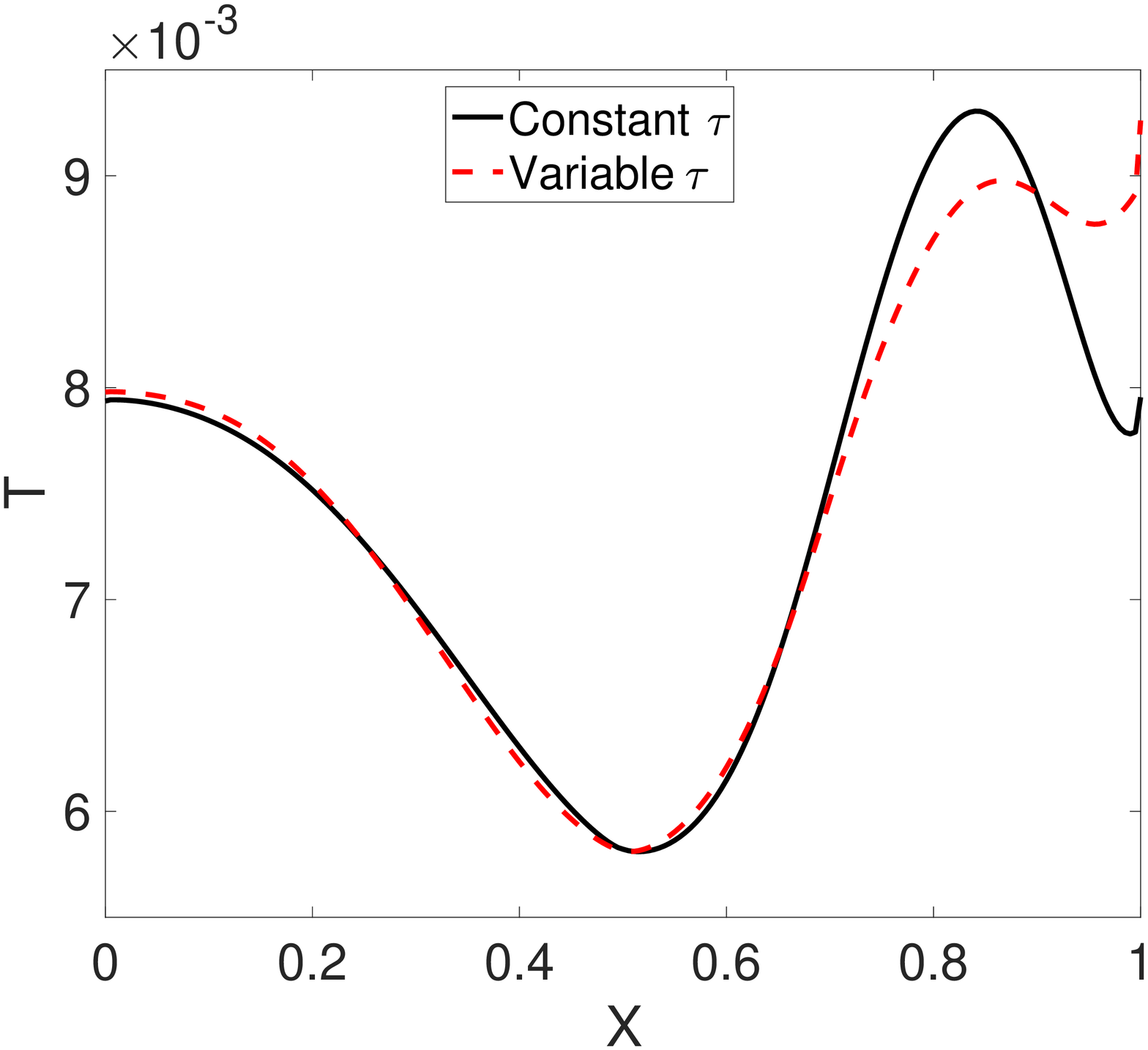}
\includegraphics[keepaspectratio=true, width=.3\textwidth]{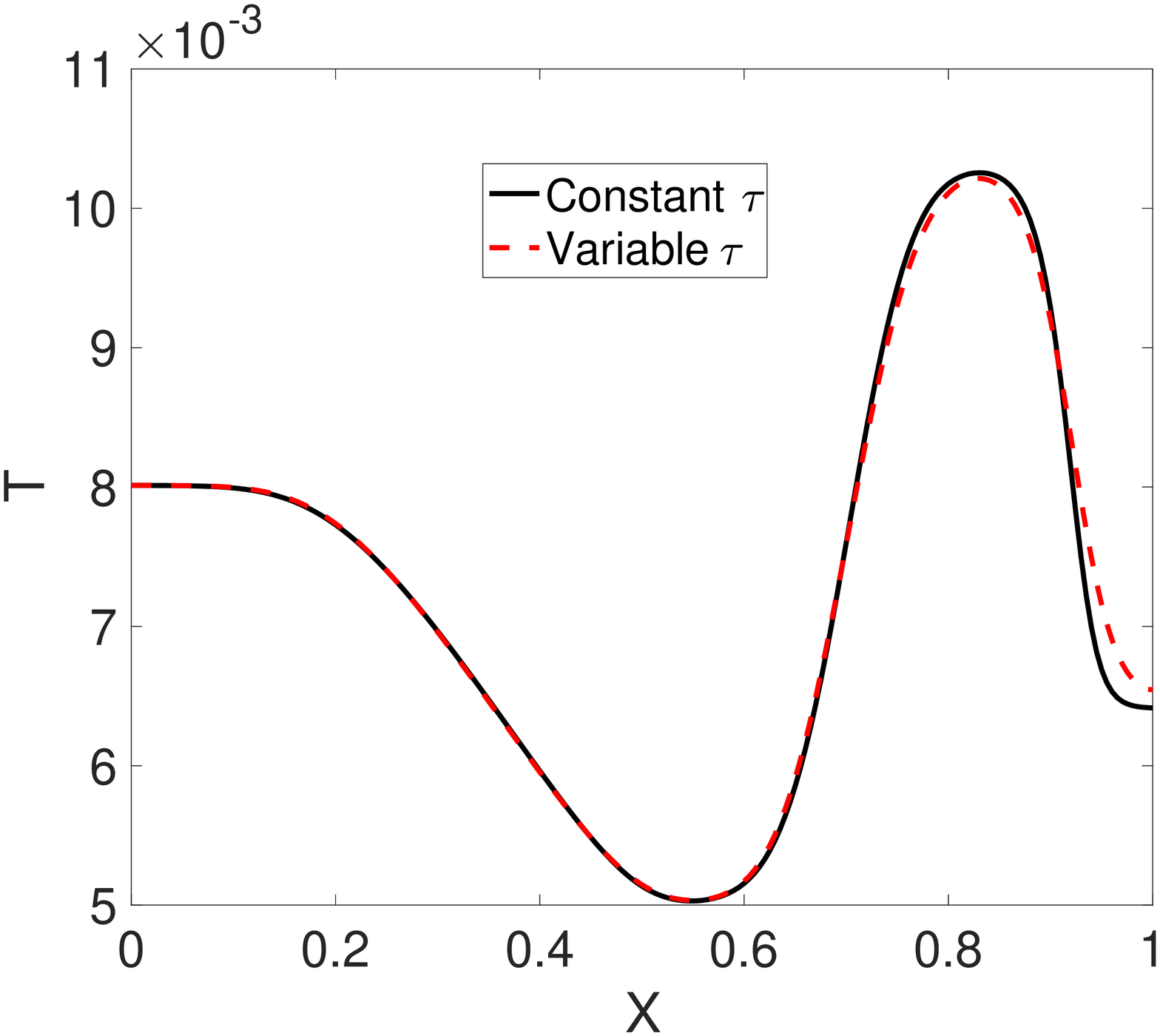}
\includegraphics[keepaspectratio=true, width=.3\textwidth]{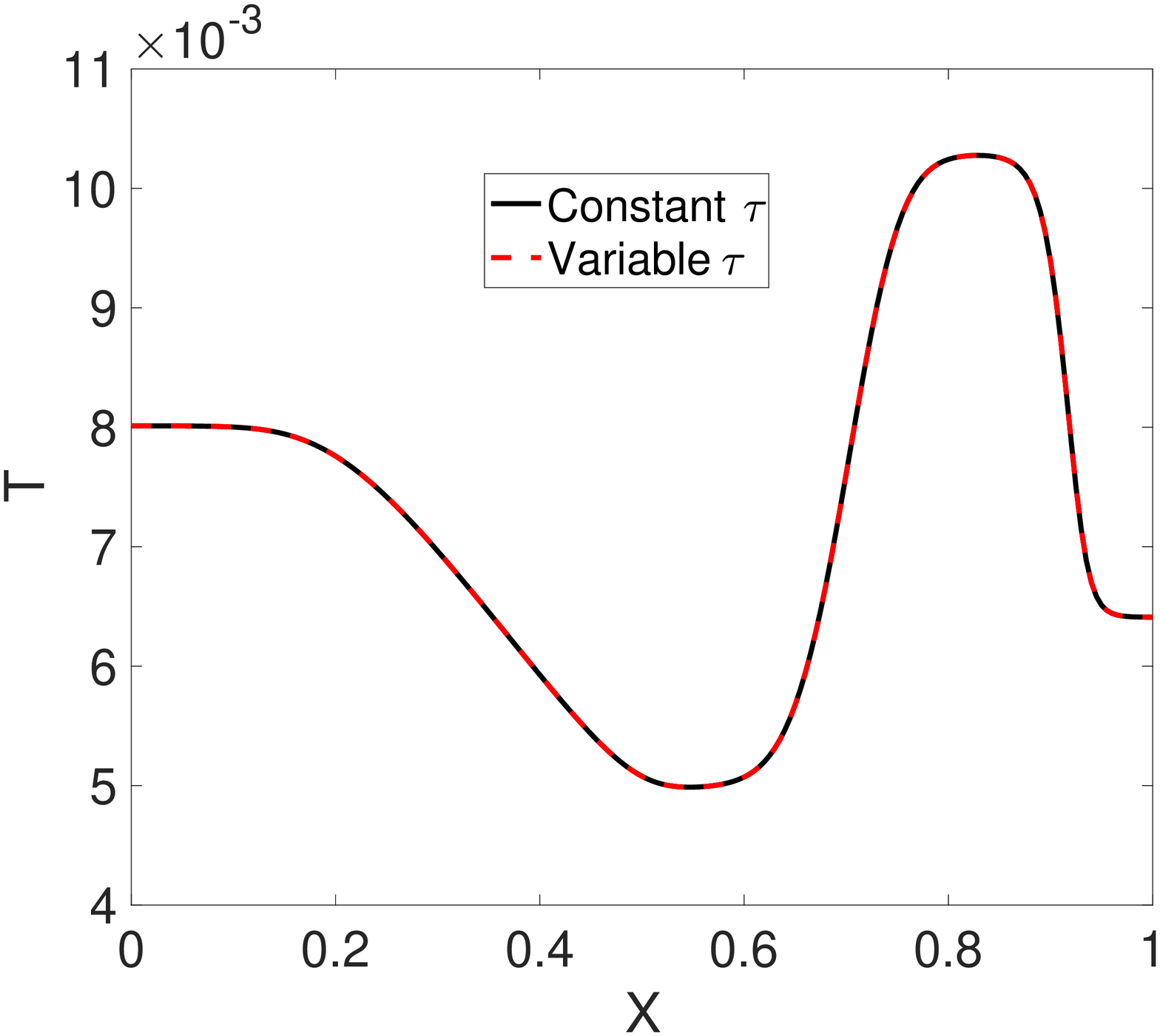}
\caption{Temperature obtained from constant and variable $\tau$ for the initial mean free paths $\lambda_l=0.02 {\rm m}, \lambda_r = 0.17$m (Left),  $\lambda_l=0.001 {\rm m}, \lambda_r = 0.008$ m  (Middle) and $\lambda_l = 10^{-7}{\rm m}$, $\lambda_r = 8\times 10^{-7}$ m (Right) with CFL = 1, $N_x = 200$ at $t_{\rm final} = 0.17$. }
\label{const_var_tau}
\end{center}
\end{figure}

In most of the DSMC simulations constant mean free paths are considered. The considerations of variable mean free paths may effect too much in the solutions since there is large fluctuation in the density in the DSMC simulations. Since we validate our 
numerical scheme for the BGK equation with the full Boltzmann equation with the help of the DSMC simulations, we fix the Knudsen number initially and keep it constant 
until the final time step. In shock tube problems there is jump in the initial mean free path. Therefore, we consider the average value of the mean free path in all time steps. 
Based on the average mean free path we define the initial relaxation time $\tau_l$ and $\tau_r$ with the help of Eq. (\ref{tau}) and use these values, for example, 
 $\tau_l$ in the domain $0\le x \le 0.5$ and $\tau_r$ in the rest of the domain throughout the simulations.  
In \cite{RF} the authors have used the constant relaxation time to solve the BGK equation. However, one can use the variable $\tau$ in the time and space, see 
Eq. (\ref{SL_step3}).  For very small $\lambda$ or $\tau$ the solutions obtained from variable and constant $\tau$ do not differ much. In Figure \ref{const_var_tau} we have plotted the temperature for constant and variable $\tau$. We have considered the density ratios $\rho_l/\rho_r = 8$. Three densities are considered, which are $\rho_l = 5\times 10^{-6} {\rm Kg}\, m^{-3}, ~ 10^{-4} {\rm Kg}\, m^{-3}, ~ 1 {\rm Kg}\, m^{-3}$ which correspond to the mean free paths $\lambda_l = 0.02, ~0.001, ~10^{-7}$ meters, respectively. We see that the difference in the solution 
obtained from the constant and variable $\tau$ for larger mean free paths. When the mean free path is very small there is no difference between the solutions.  

\subsection{Test 3: Comparison of solutions in regular vs irregular grids}

\begin{figure}
\includegraphics[keepaspectratio=true, width=.48\textwidth]{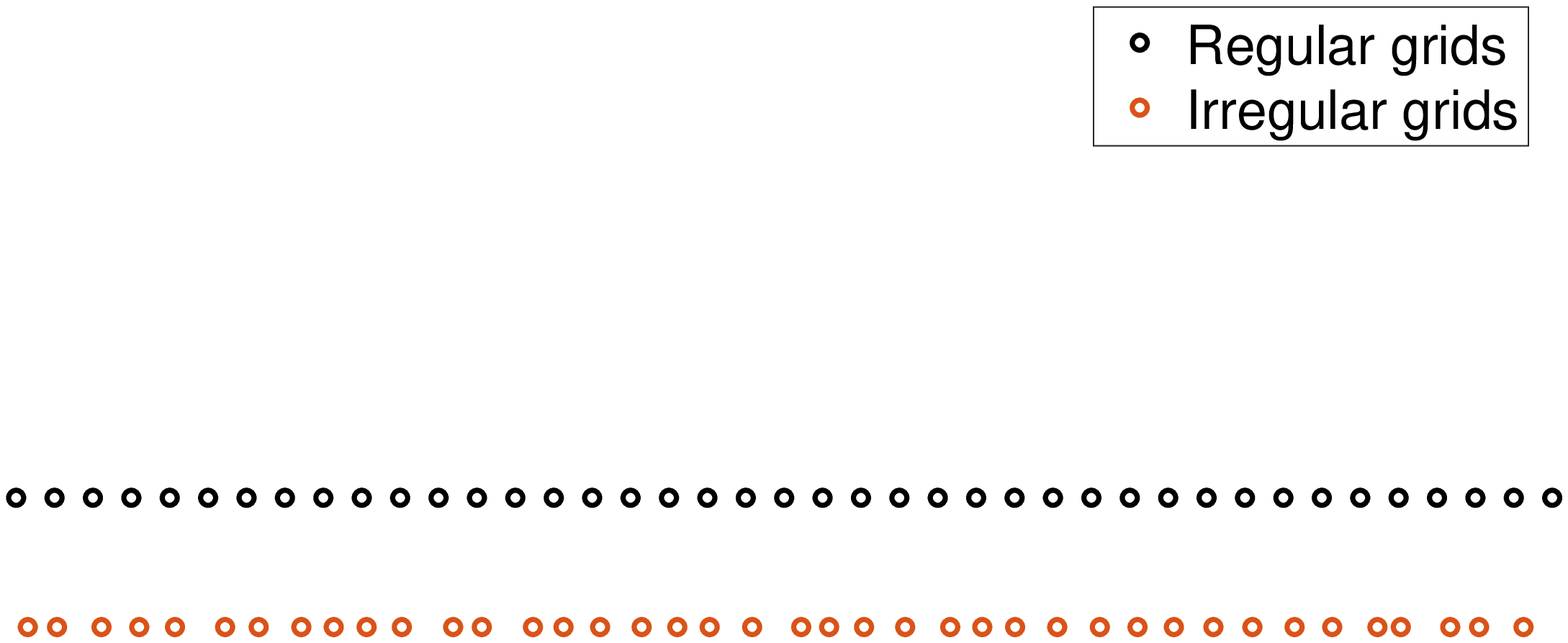}
\includegraphics[keepaspectratio=true, width=.48\textwidth]{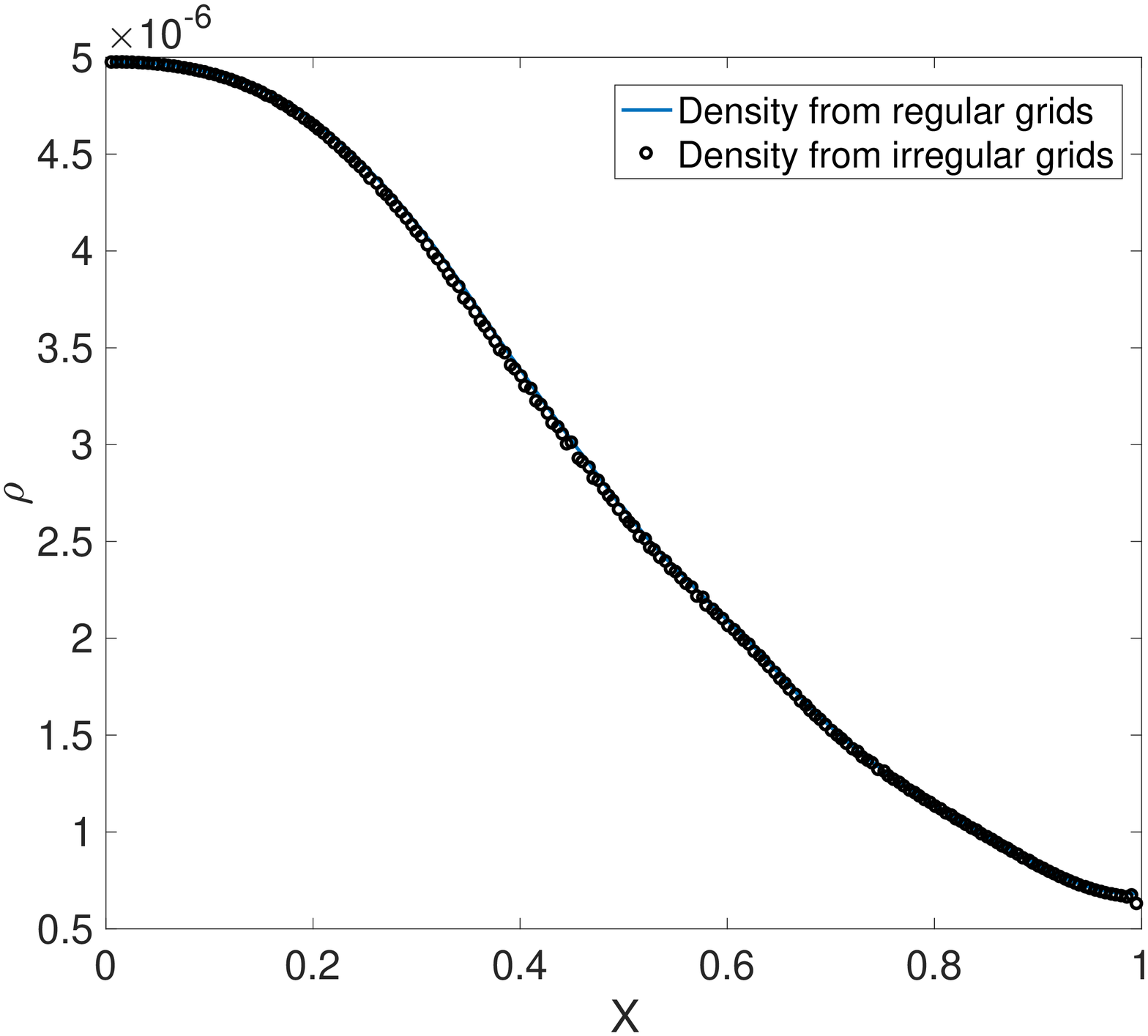}
\caption{Zoom of regular and irregular grids (left) and comparison of density for regular vs irregular grids for $N_x = 200$ CFL = 1, initial mean free path $\lambda_l = 0.02$ m for $x<0.5$,  $\lambda_r = 0.17$ m for $x > 0.5$ and $t_{\rm final} = 0.17$. }
\label{reg_vs_irre}
\end{figure}

In this test case we present the comparison of the numerical solutions in regular as well as irregular grids. We have considered $N_x = 200$.  
The regular grids are generated according to $x_i= (i-1)* \Delta x, i = 1, \ldots, N_x+1$. To create the irregular grids we have moved the  
regular grids with velocity $\Delta x/4$ times the random number  $i=2, \ldots, N_x$. This movement is performed for 2 iterations. 
The densities are 
$\rho_l = 5\times 10^{-6} {\rm Kg}~m^{-3}$ and $\rho_r = 0.125\times 10^{-6} {\rm Kg}~m^{-3}$. The corresponding initial mean free paths are $\lambda_l = 0.02$ m on the left and $\lambda_r = 0.17$ m on the right half of the domain. The corresponding initial relaxation times are $0.01$ on the left half of the domain and $0.0957$ on the right half of the domain. 
The flow is in transition regime.  In Figure \ref{reg_vs_irre} we have plotted the zoom of the regular and irregular grids on the left and 
the densities obtained from the MLS interpolation in regular as well as irregular grids. 
We observe that the irregular grids make no difference to the solutions obtained from the regular grids.

\subsection{Test 4}

In this case, we have compared the numerical solutions obtained from the BGK model with the DSMC simulations \cite{Babovsky, NS} for the Boltzmann equation since DSMC results are widely used as Benchmark solutions. 
 All parameters are same as in the Test 2. For DSMC simulations we have considered  $200$ cells and $400$ gas molecules are initially generated per cell according to the Maxwellian 
 distribution in the velocity, where the initial density, temperature and velocity are its parameters. Notice that DSMC is a method to approximate the solution to the Boltzmann equation for hard spheres, while BGK is a simplified model of the Boltzmann equation, therefore we do not expect to observe the same behaviour. Both models have in common that in the 
limit of very small mean free path converge to the compressible Euler equations for a monoatomic gas, therefore we expect that the two models 
provide similar results when adopted with the same small mean free path. Furthermore, we observe that standard BGK fails to correctly capture first order effects in the (small) mean free path, because it is not possible to match at the same time thermal conductivity and viscosity coefficient with the single parameter $\tau$. There are extensions of the BGK model, such as, for example, the so called 
ES-BGK.  It is possible to mention the existence of other BGK models, such as for example the ellipsoidal BGK \cite{Ho},  that are able to capture the correct Navier-Stokes limit. The use of such models is however beyond the scope of the present paper. 
 
\begin{figure}
\includegraphics[keepaspectratio=true, width=.329\textwidth]{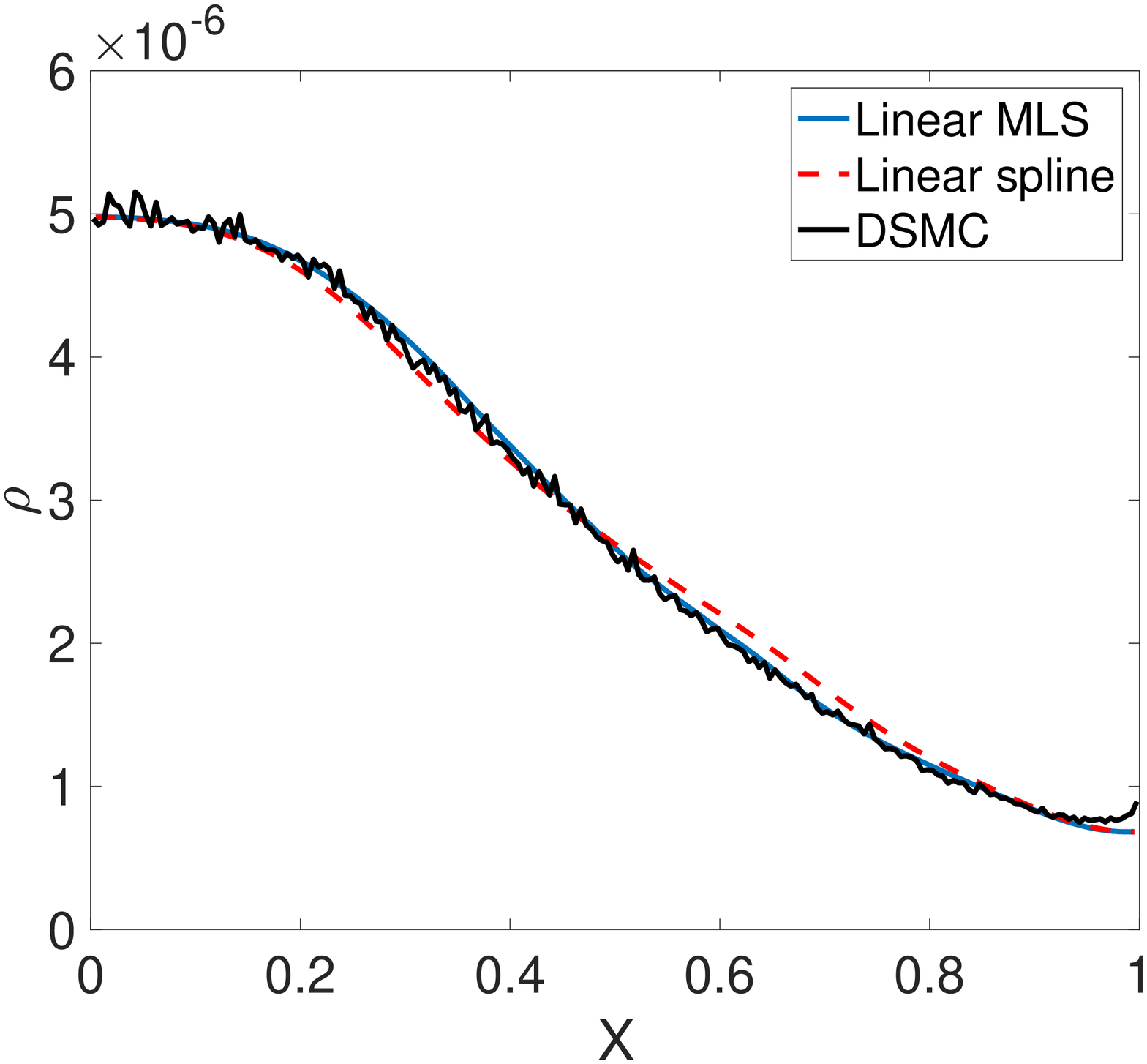}
\includegraphics[keepaspectratio=true, width=.329\textwidth]{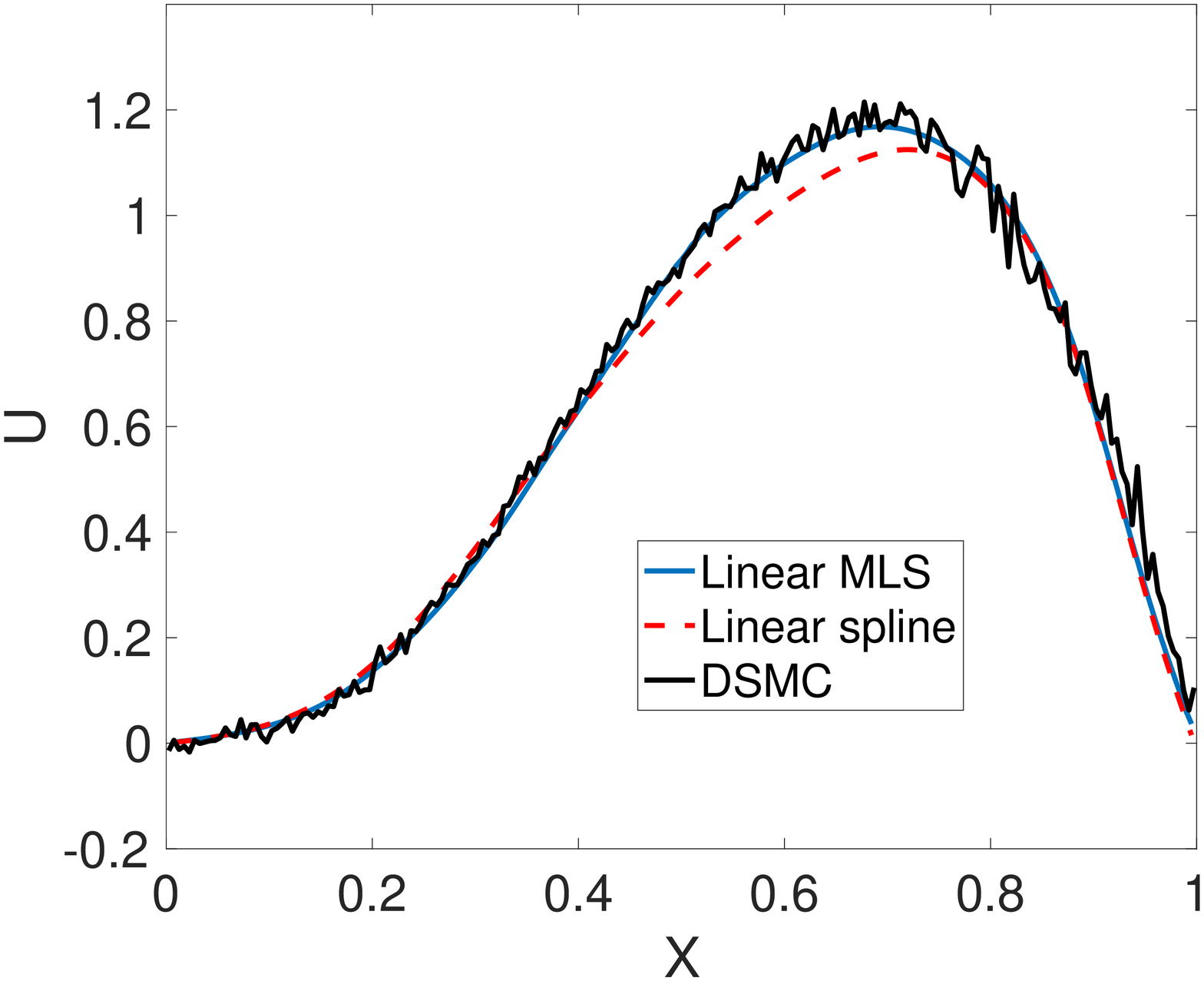}
\includegraphics[keepaspectratio=true, width=.329\textwidth]{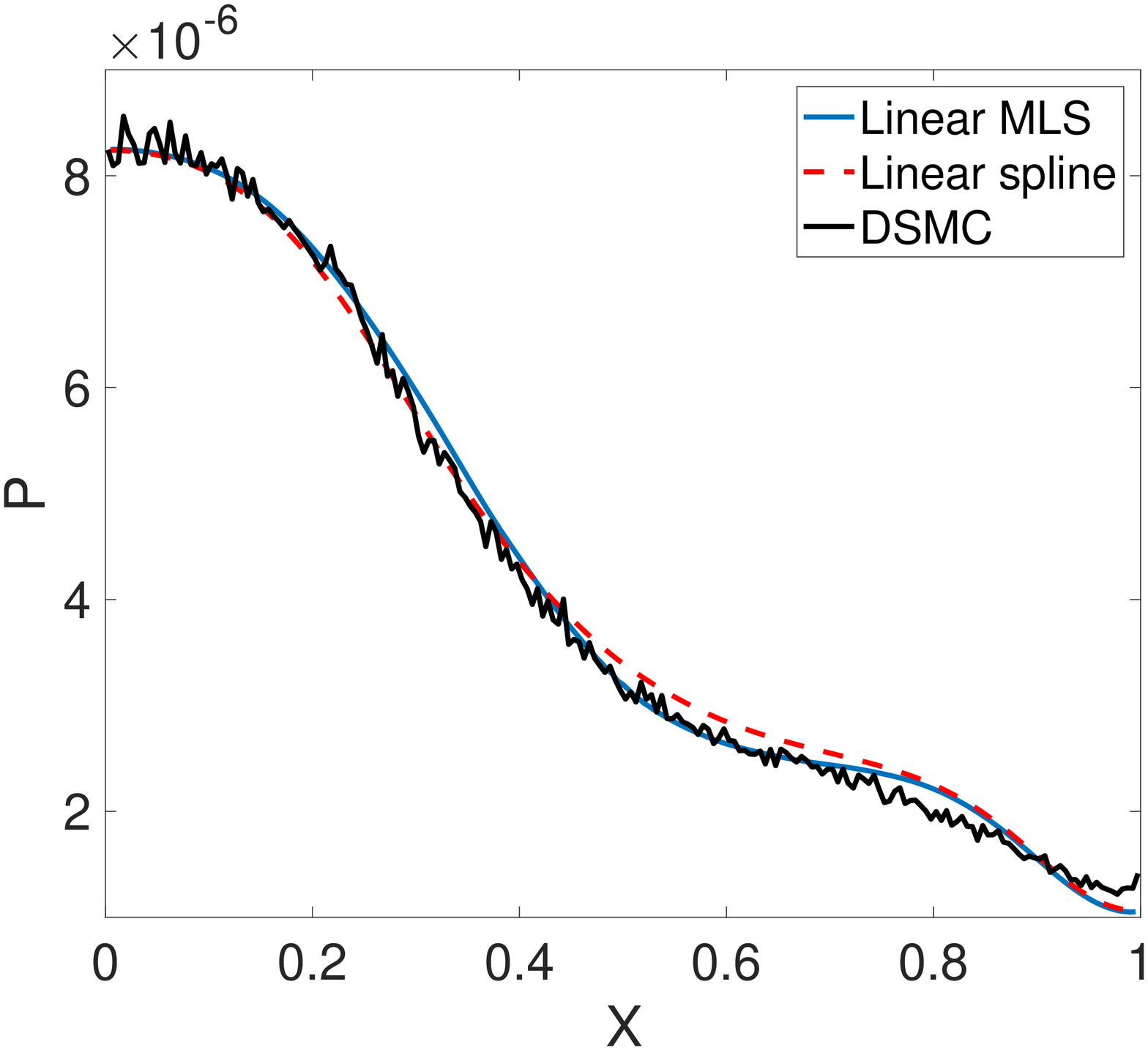}
\caption{Comparison of density, velocity and pressure obtained from MLS reconstruction, Spline interpolation and DSMC for initial mean free path is $\lambda_l = 0.02$ m for $x<0.5$ and $\lambda_r = 0.17$ m for $x>0.5$ for CFL = 1, $Nx = 200$, and $t_{\rm final} = 0.17$.}
\label{test2}
\end{figure}

 The cell size for DSMC simulations must be smaller than the mean free path. Therefore, the time step of the DSMC solver is restricted by the mean free path, see \cite{Babovsky, NS} for details. However, the time step for the Semi Lagrangian scheme is not restricted, which has shown in Test 1. 
 Since the flow has low Mach number, the statistical fluctuations dominates the DSMC solutions. 
 Therefore, we have obtained $20$ independent runs. 
 The initial density is chosen same as in the Test 3. 
In Figure \ref{test2} we have plotted the density, velocity and pressure obtained from the DSMC simulations and the BGK model using linear MLS and linear spline interpolations. We observe that the linear MLS interpolation scheme 
for the BGK model and DSMC results match perfectly, however,  the linear spline interpolation gives some deviations from the solutions of the DSMC simulations. 

\subsection{Test 5}

In this test case we have increased the density by factor $2$ and $20$ times compared to the Test case 3 and 4 such that $\rho_l =  10^{-5} {\rm Kg}m^{-3}$ and $ 10^{-4} {\rm Kg}m^{-3}$, respectively. 
The corresponding mean free paths are $0.01$ and $0.001$ meters on the left half of the domain and $8$ times larger on the right half.  For $\rho_l =10^{-5} {\rm Kg} m^{-3}$ we have considered the $400$ cells in the DSMC simulations and the same number 
of grids for the BGK model. Similarly, for $\rho_l =10^{-4} {\rm Kg} m^{-3}$ we have considered the $1000$ cells in the DSMC simulations and the  $800$ grids for the BGK model. The increase of the number of cells 
is due to the restriction that the DSMC cells must be smaller than the mean free path. Other parameters are same as in the earlier test cases. 

\begin{figure}
\includegraphics[keepaspectratio=true, width=.329\textwidth]{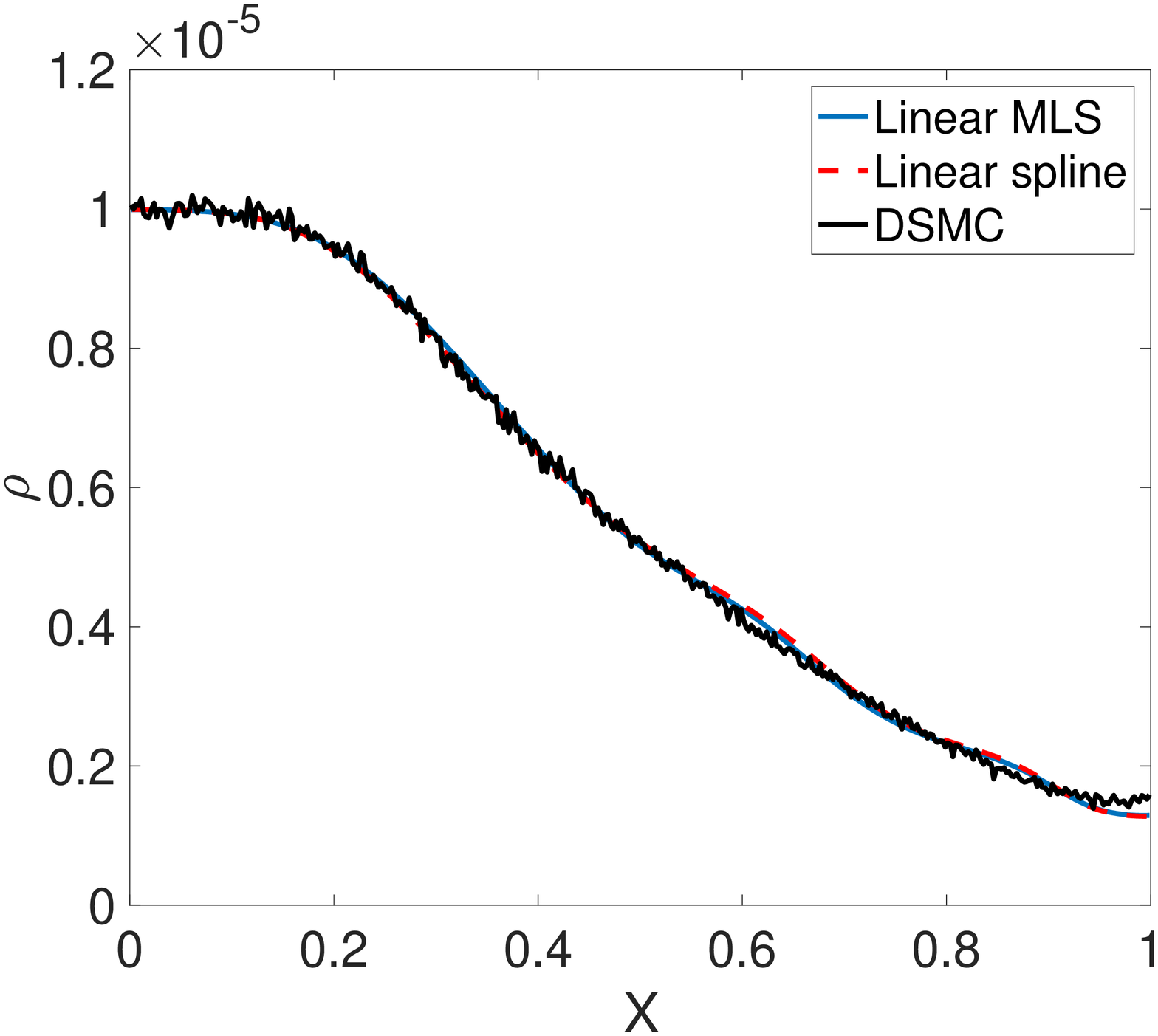}
\includegraphics[keepaspectratio=true, width=.329\textwidth]{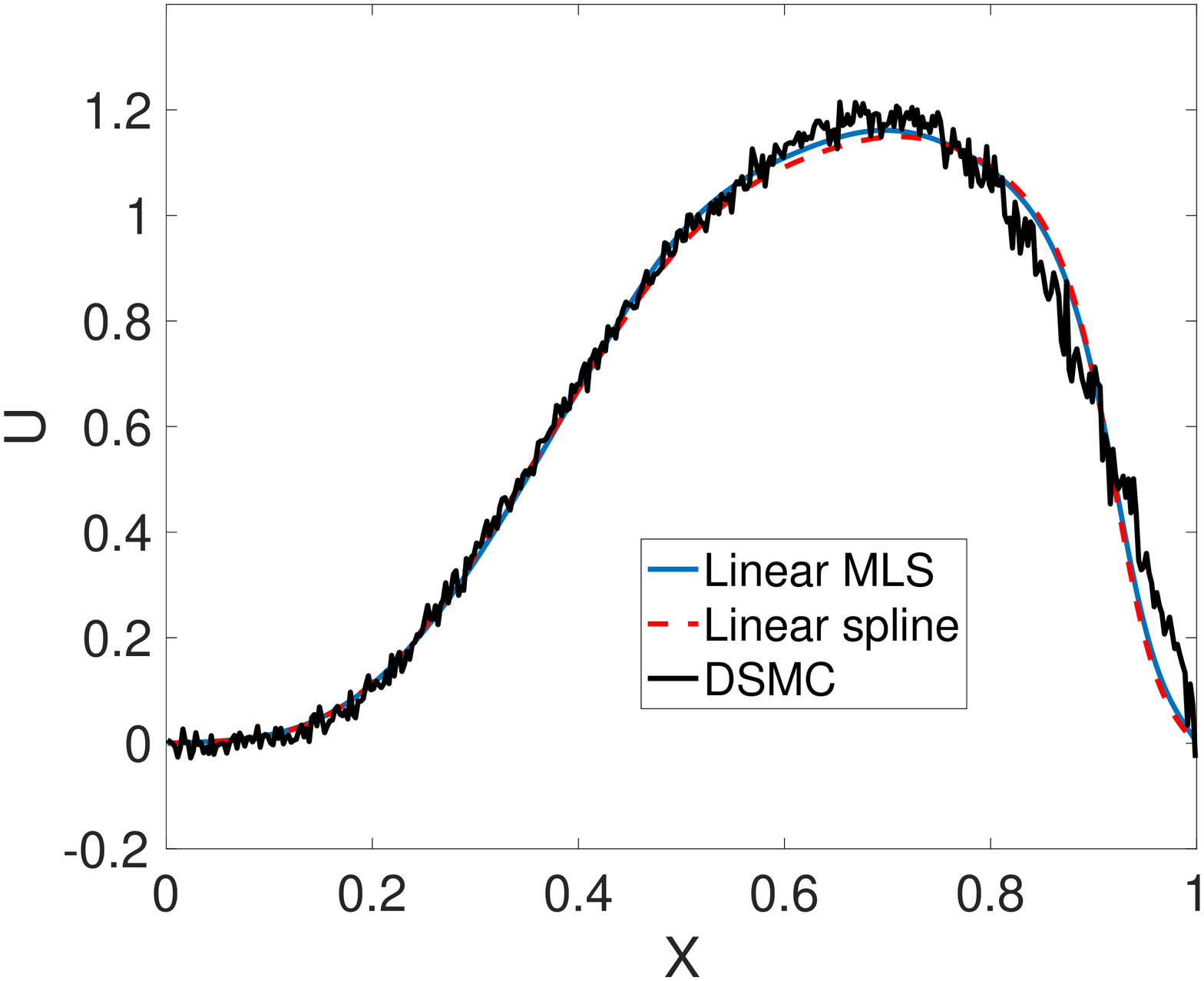}
\includegraphics[keepaspectratio=true, width=.329\textwidth]{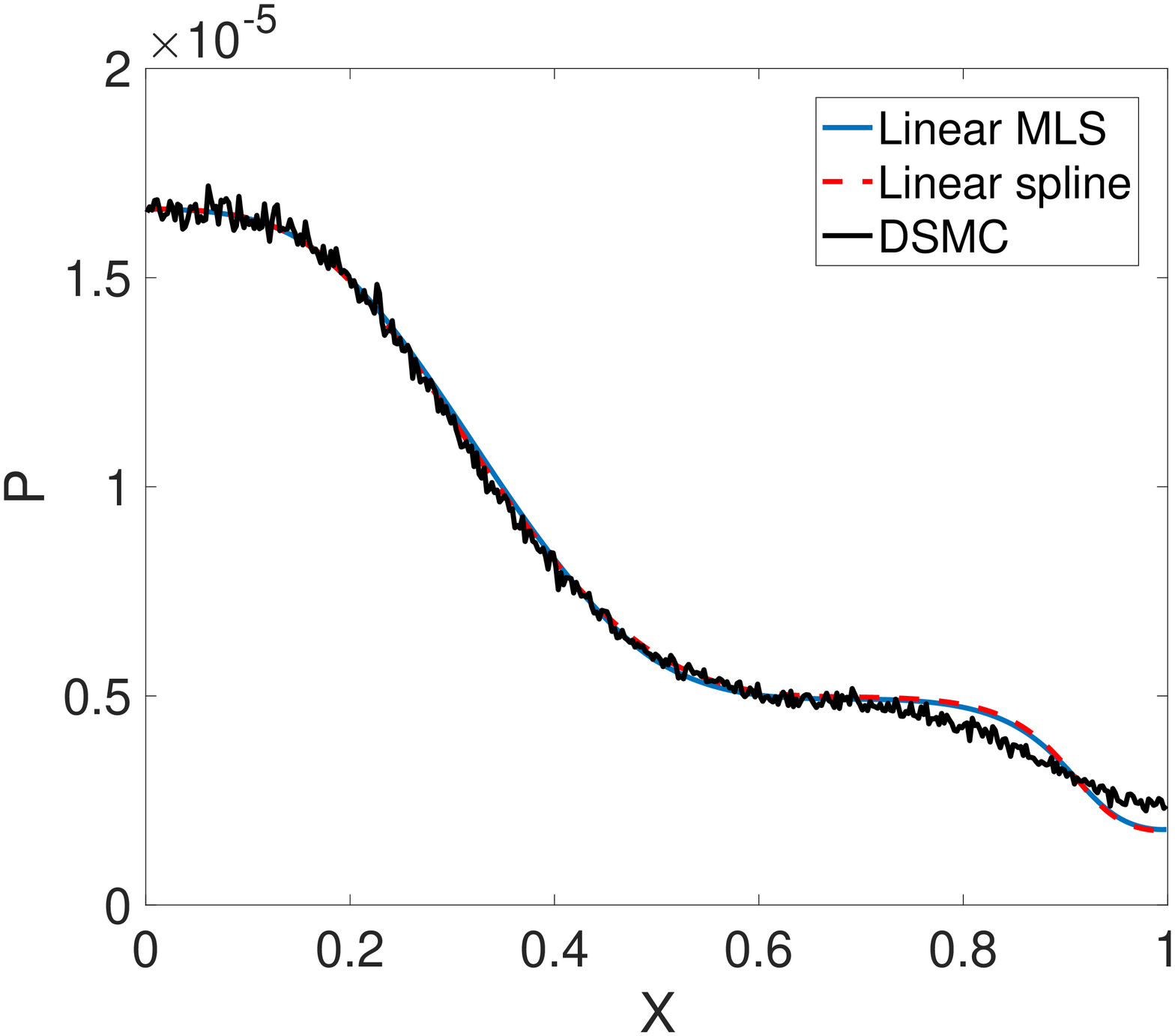}
\caption{Comparison of density, velocity and pressure obtained from MLS reconstruction, Spline interpolation and DSMC for initial mean free path $\lambda_l = 0.01$ m for $x < 0.5$,  $\lambda_r = 0.08$ m for $x>0.5$ for CFL = 1, $Nx = 400$, $t_{\rm final} = 0.17$.}
\label{test3a}
\end{figure}

\begin{figure}
\includegraphics[keepaspectratio=true, width=.329\textwidth]{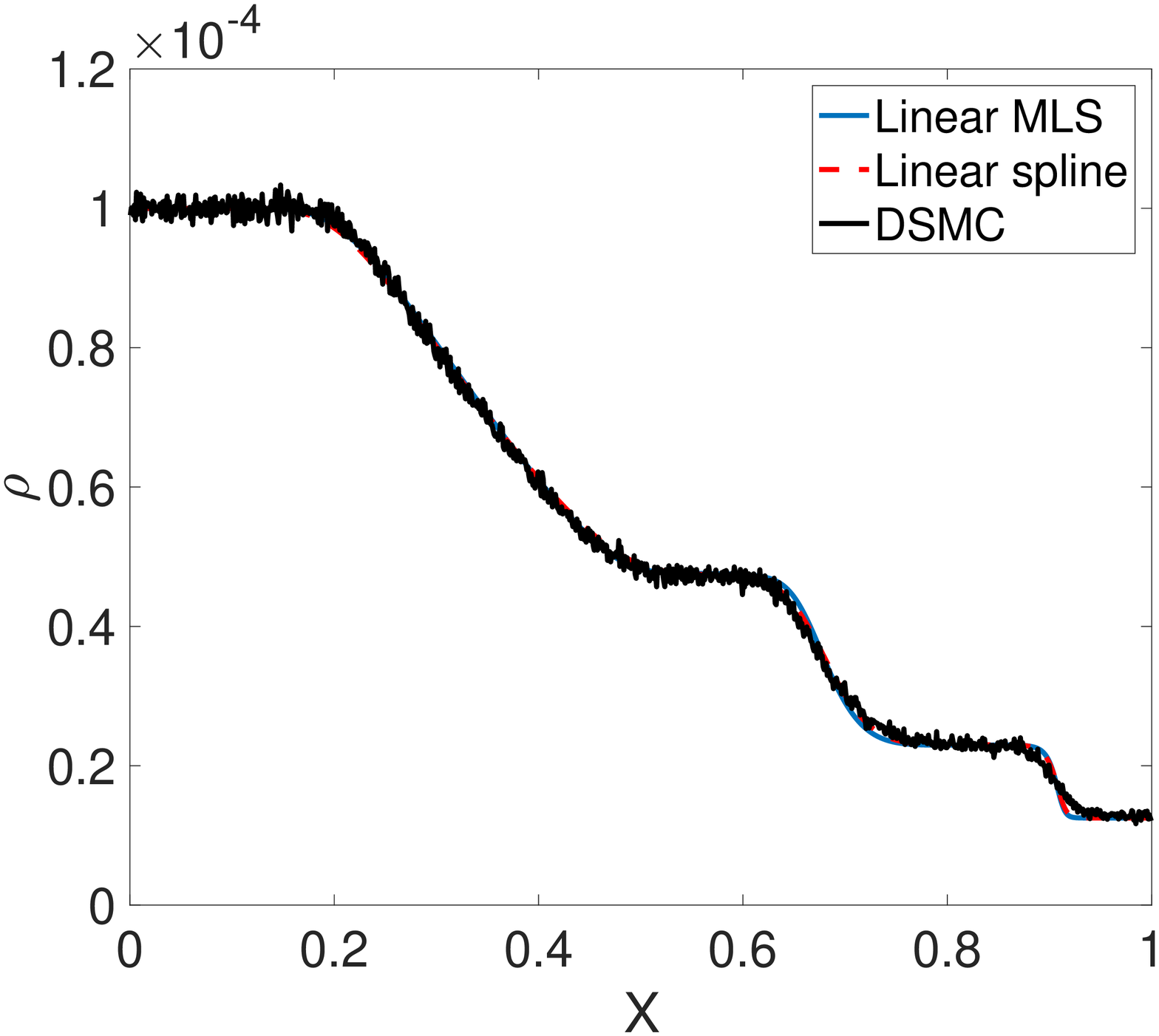}
\includegraphics[keepaspectratio=true, width=.329\textwidth]{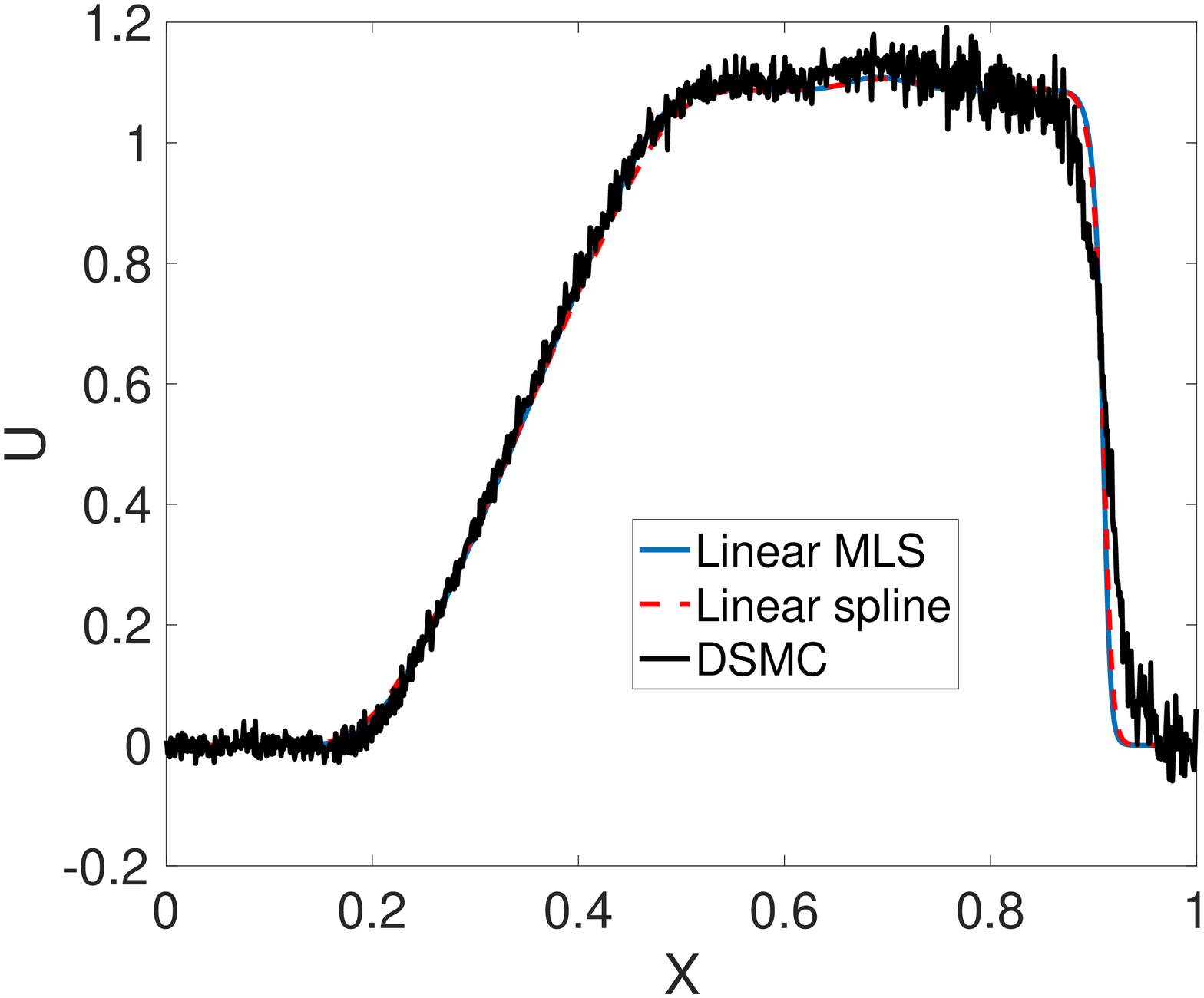}
\includegraphics[keepaspectratio=true, width=.329\textwidth]{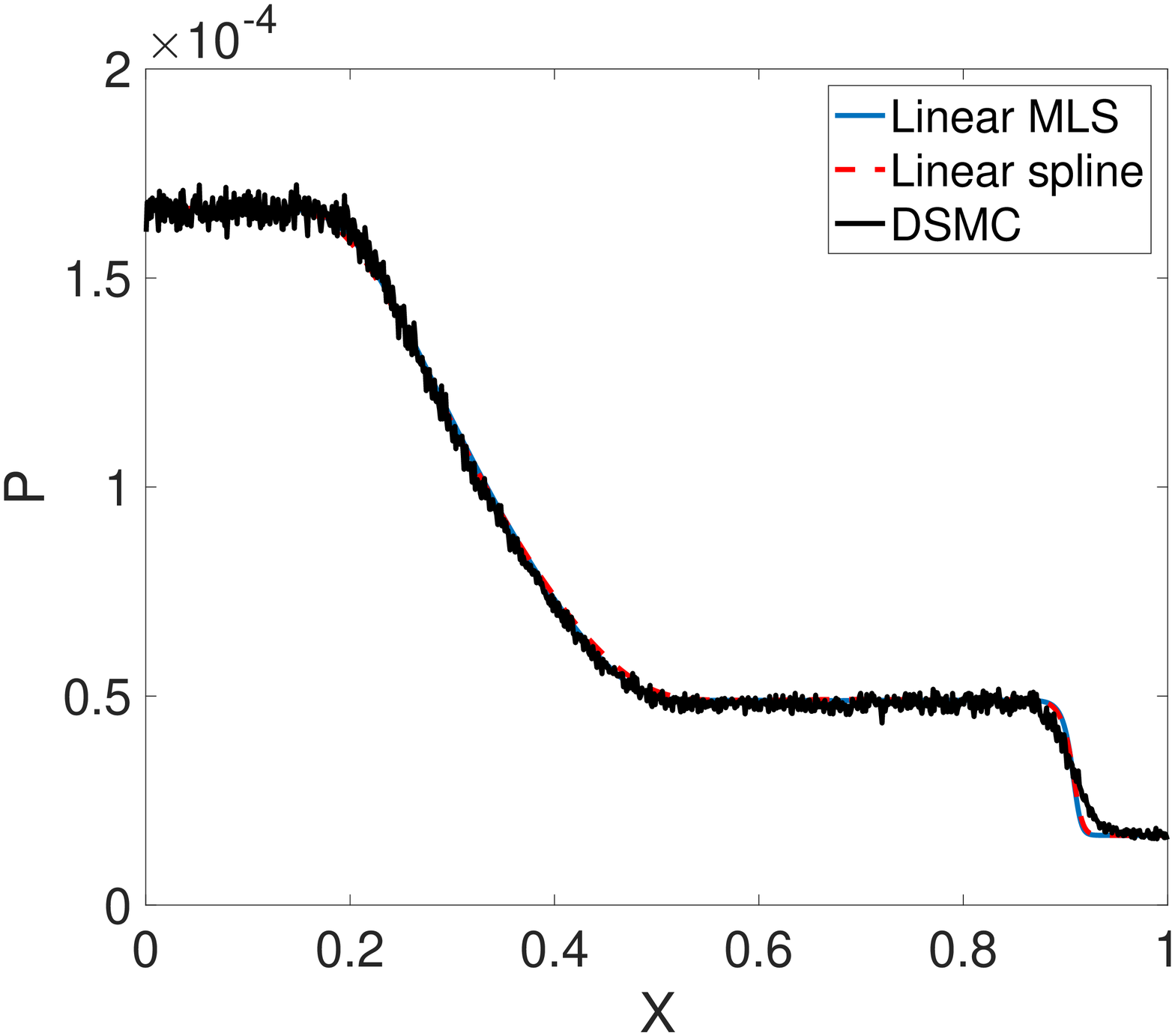}
\caption{Comparison of density, velocity and pressure obtained from MLS reconstruction, Spline interpolation and DSMC for initial mean free path $\lambda_l = 0.001$ m for $x<0.5$ and $\lambda_r = 0.008$ m for $x > 0.5$, $Nx = 800$ for the BGK and $Nx=1000$ for the DSMC, CFL = 1, $t_{\rm final} = 0.17$. }
\label{test3b}
\end{figure}

For the DSMC simulations we have performed $20$ independent runs. In Figures \ref{test3a} and \ref{test3b} we have plotted the numerical solutions from all three methods. 
We again observe that the DSMC solutions and the solutions of the BGK model obtained by linear MLS interpolation are closer than the linear spline interpolation. 

\subsection{Test 6}

In this test case we consider  $\rho_l = 1 {\rm Kg} m^{-3}$ and $\rho_r = 0.125 {\rm Kg} m^{-3}$. The corresponding left and right mean free paths are $10^ {-7}$ m and $8\times 10^{-7}$ m, respectively. 
We note that the size of DSMC cells must be smaller than $\lambda$, so, we need at least $9\times 10^{6}$ cells for such mean free path. The time step also has to be reduced accordingly. Moreover, 
the number of gas molecules is also very high and the computational time for the DSMC simulations becomes enormously high. 
In one dimensional case, this it quite a large number of cells. Therefore, we have not performed DSMC simulations in this case. However, there is no restriction of cell size for the Semi Lagrangian scheme for 
the BGK model. We have again  used $800$ grids for the Semi Lagrangian scheme for this smaller mean free paths. 

\begin{figure}
\includegraphics[keepaspectratio=true, width=.329\textwidth]{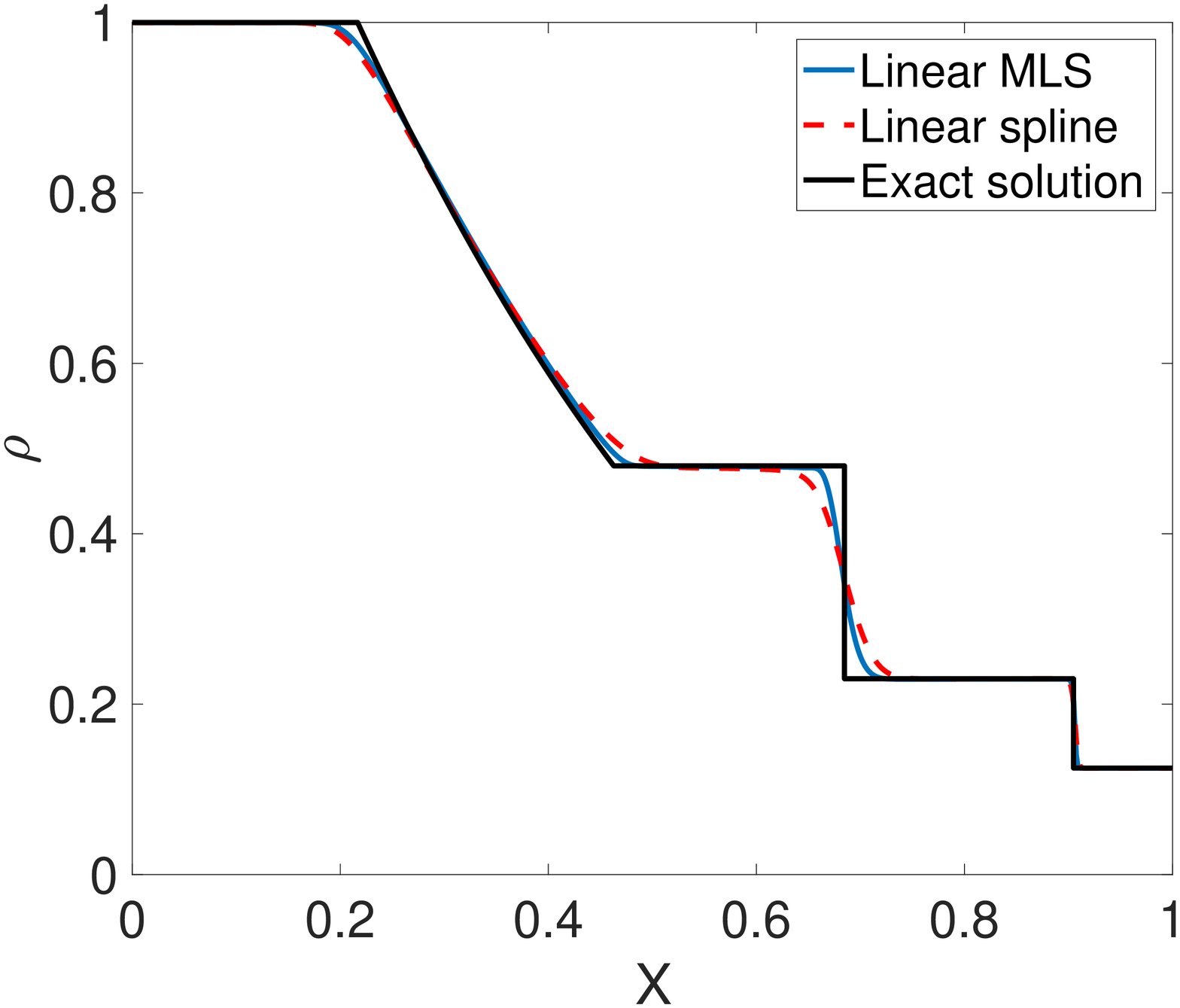}
\includegraphics[keepaspectratio=true, width=.329\textwidth]{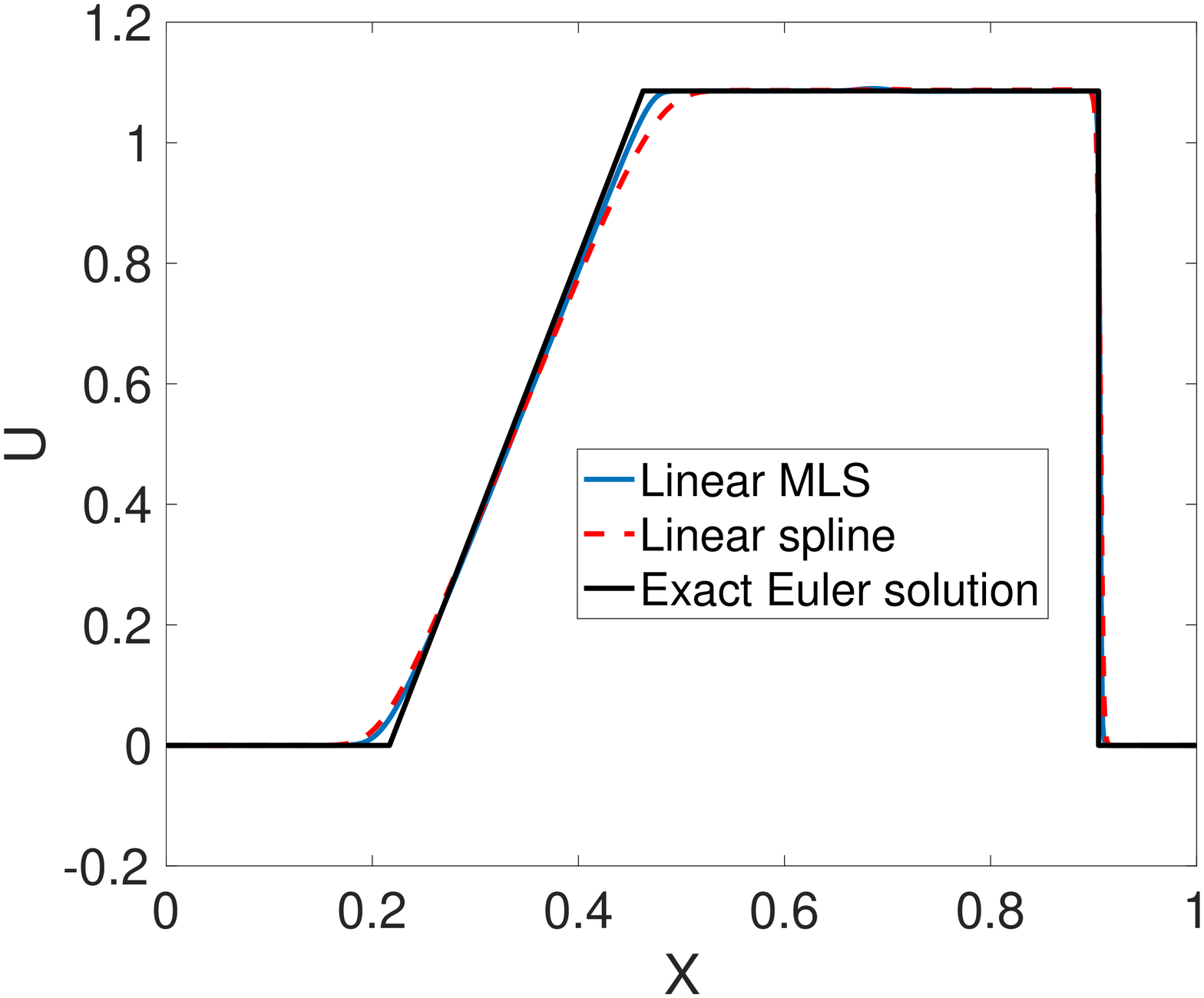}
\includegraphics[keepaspectratio=true, width=.329\textwidth]{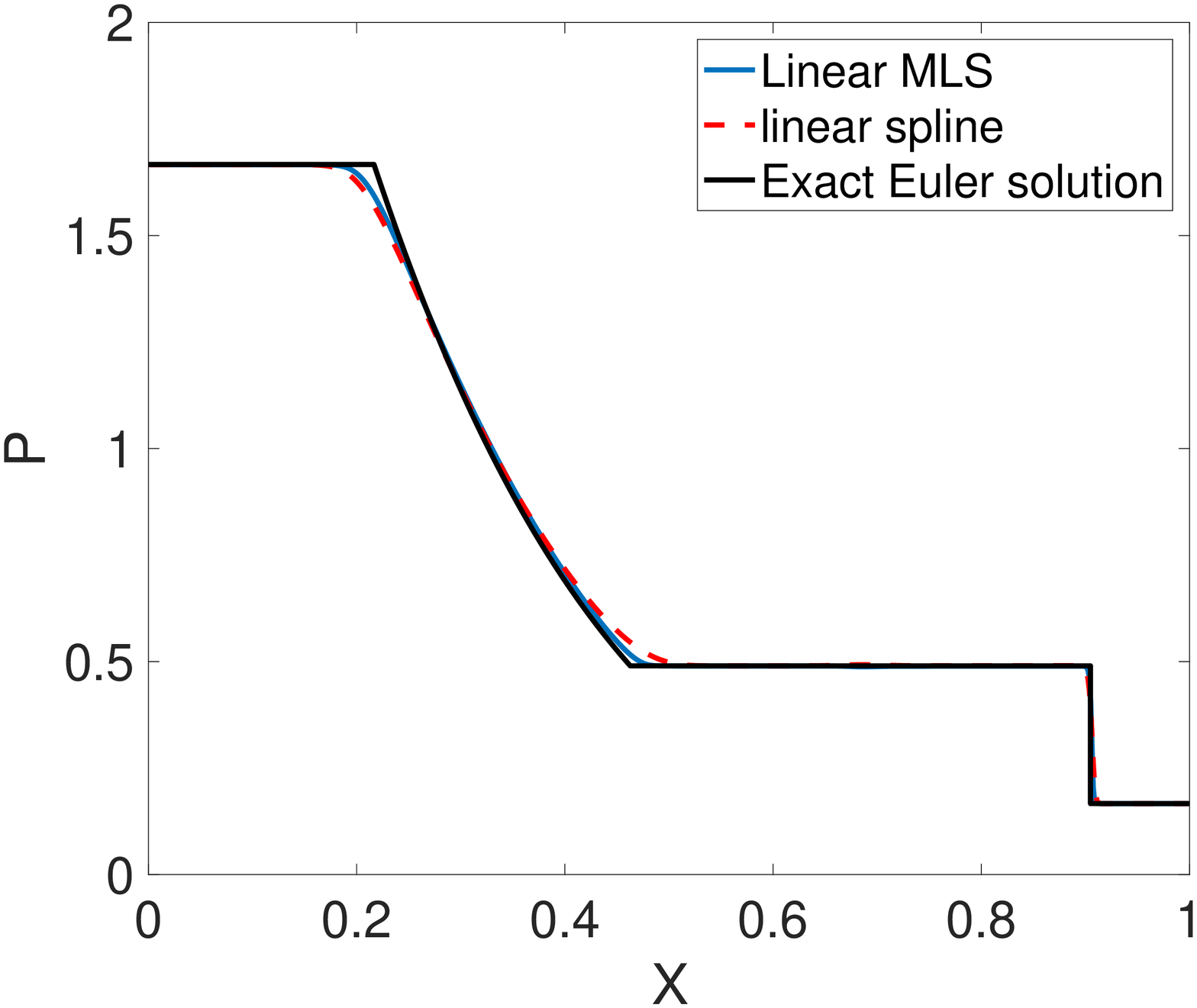}
\caption{Comparison of density, velocity and pressure obtained from MLS reconstruction, Spline interpolation and the exact solutions of the compressible equations for initial mean free path $\lambda_l = 10^{-7}$ m for $x<0.5$ and $\lambda_r = 8\times 10^{-7}$ m for $x > 0.5$, $Nx = 800$, CFL = 1, $t_{\rm final} = 0.17$. }
\label{compare_euler}
\end{figure}

On the other hand for this small mean free path we can solve the continuum equations, for example,  the compressible Euler equations. For 
the shock tube problem, the compressible Euler equations can be solved exactly. In Figure \ref{compare_euler} we have again plotted the density, velocity and pressure at final time $0.17$ seconds 
obtained from the MLS and Spline interpolations together with the exact solutions of the compressible Euler equations. In the shock region all three solutions match perfectly, however, in the contact discontinuity and the rarefaction region, the 
linear MLS interpolation scheme gives better approximation than the linear spline interpolation scheme. 

 \subsection{Test 7: Comparing solutions from continuous and discrete Maxwellian}

In all above examples we have considered the velocity grids $N_v = 20$  in all three directions the standard Maxwellian given by (\ref{maxwellian}). Increasing the number of velocity grids do not help much in the accuracy of the solutions. However, smaller 
values of $N_v$ affects the accuracy. This means we loose the conservative properties of the scheme. 
In this test case we have considered all parameters as in the Test 6. We have observed that $N_v = 20$ gives solutions close as the exact solutions of the compressible 
Euler equations. So, the solutions obtained from $N_v = 20$ are our reference solutions. 
We have  decreased the values of $N_v$ and compare the solutions with the reference solutions.  The smallest one which gives the stable solutions is $N_v = 13$ for the case of the 
standard Maxwellian. 
But the solutions deviate from the reference solutions. For smaller values of $N_v$ we loose the conservative properties. In other words,  the moments obtained from the discrete summation in 
Eq. (\ref{discrete_moments1} - \ref{discrete_moments5}) are not exactly equal to the moments computed from the standard Maxwellian (\ref{maxwellian}). 
To obtain conservative properties, one uses the discrete Maxwellian suggested by Mieussiens \cite{Mieuss}. The discrete Maxwellian depends on five parameters, which are also related to the moments and the parameters are determined by solving the non-linear system of equations. Since the Jacobian matrix has very bad condition number, the standard Newton's method for solving non-linear system does not work, one has to use the back tracking line search algorithm \cite{DS}. 
A careful comparison between the use of continuous and discrete Maxwellian in Semi Lagrangian schemes for the computation of shock problems is reported in \cite{BCRY}.

\begin{figure}
\includegraphics[keepaspectratio=true, width=.329\textwidth]{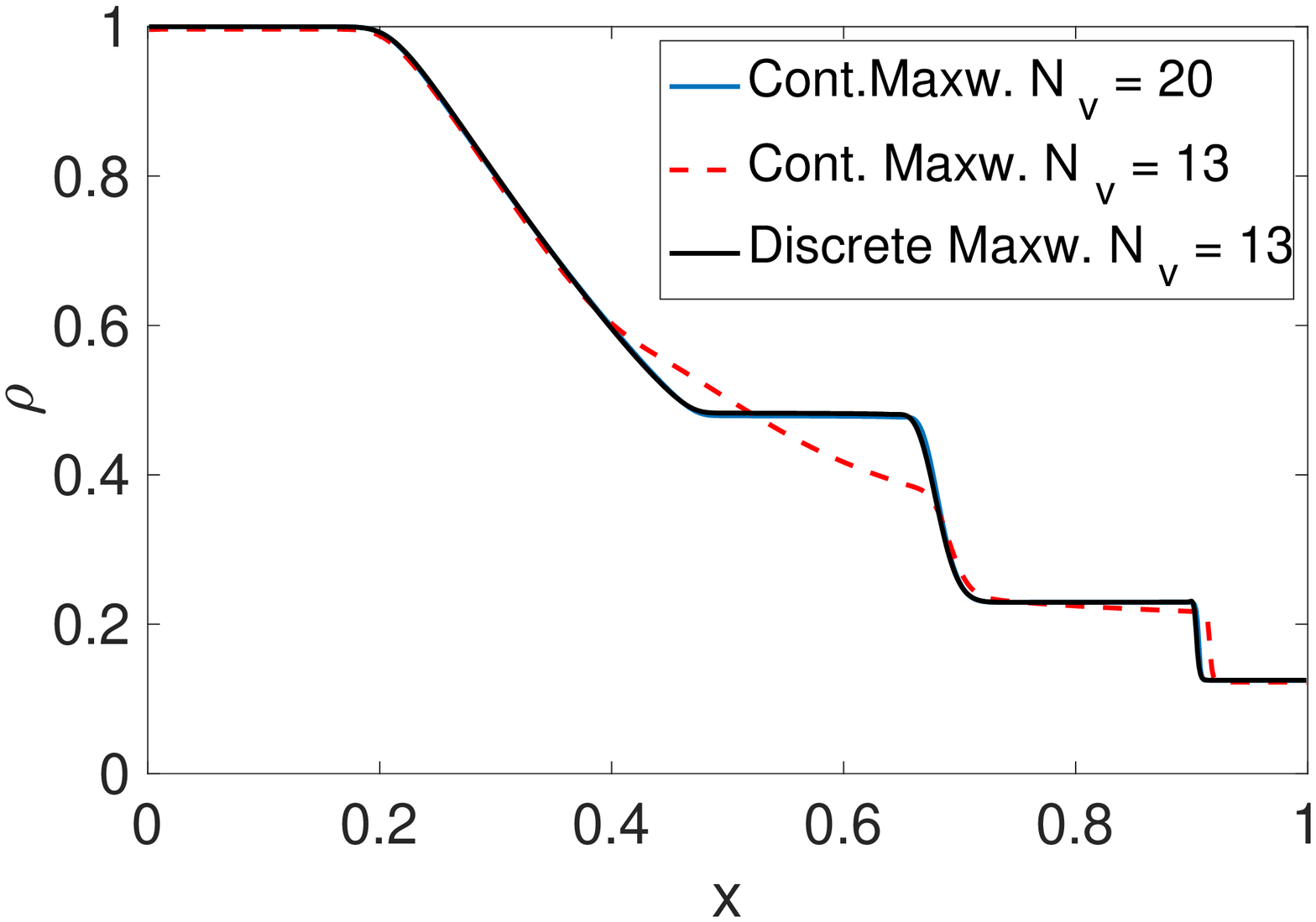}
\includegraphics[keepaspectratio=true, width=.329\textwidth]{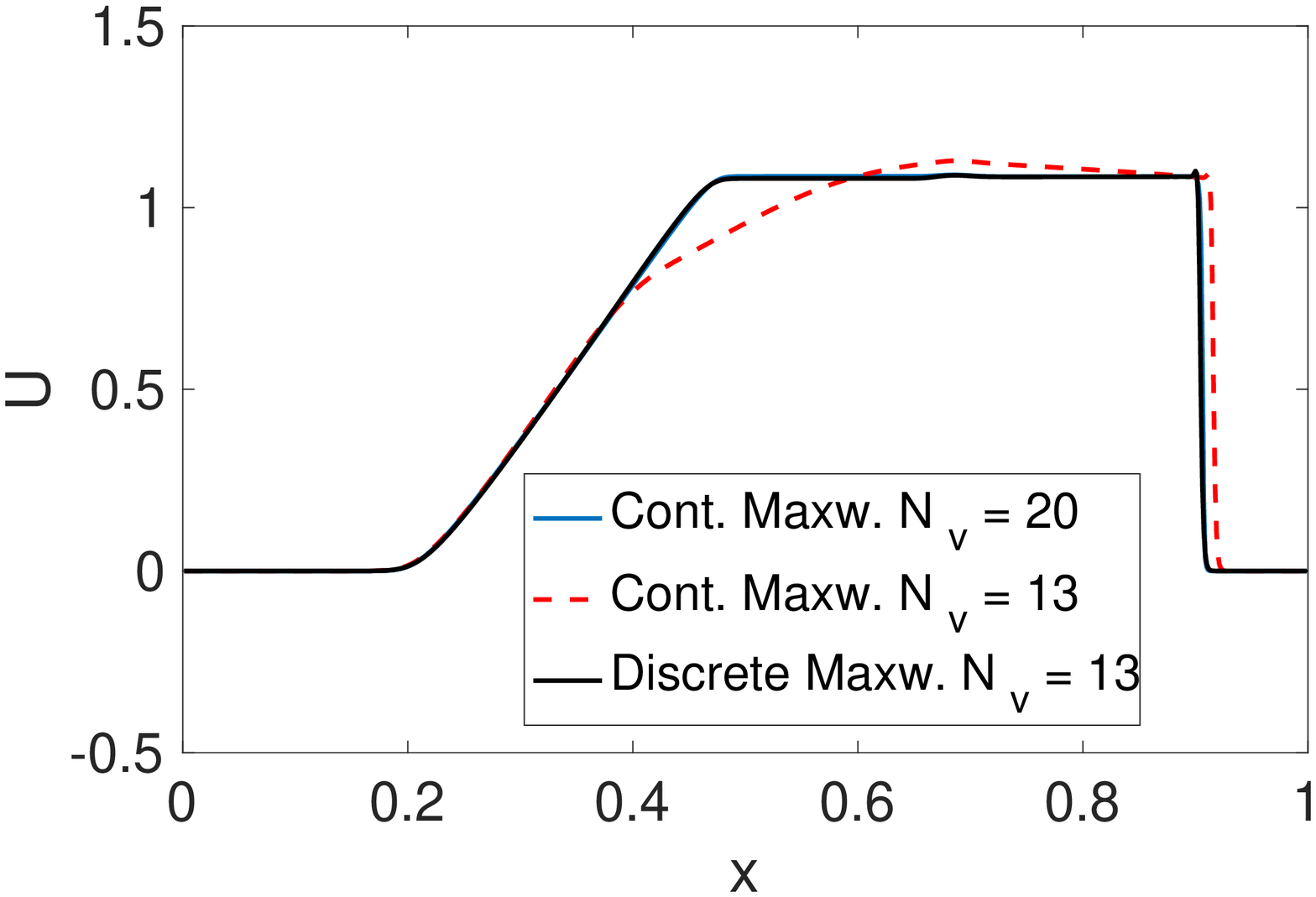}
\includegraphics[keepaspectratio=true, width=.329\textwidth]{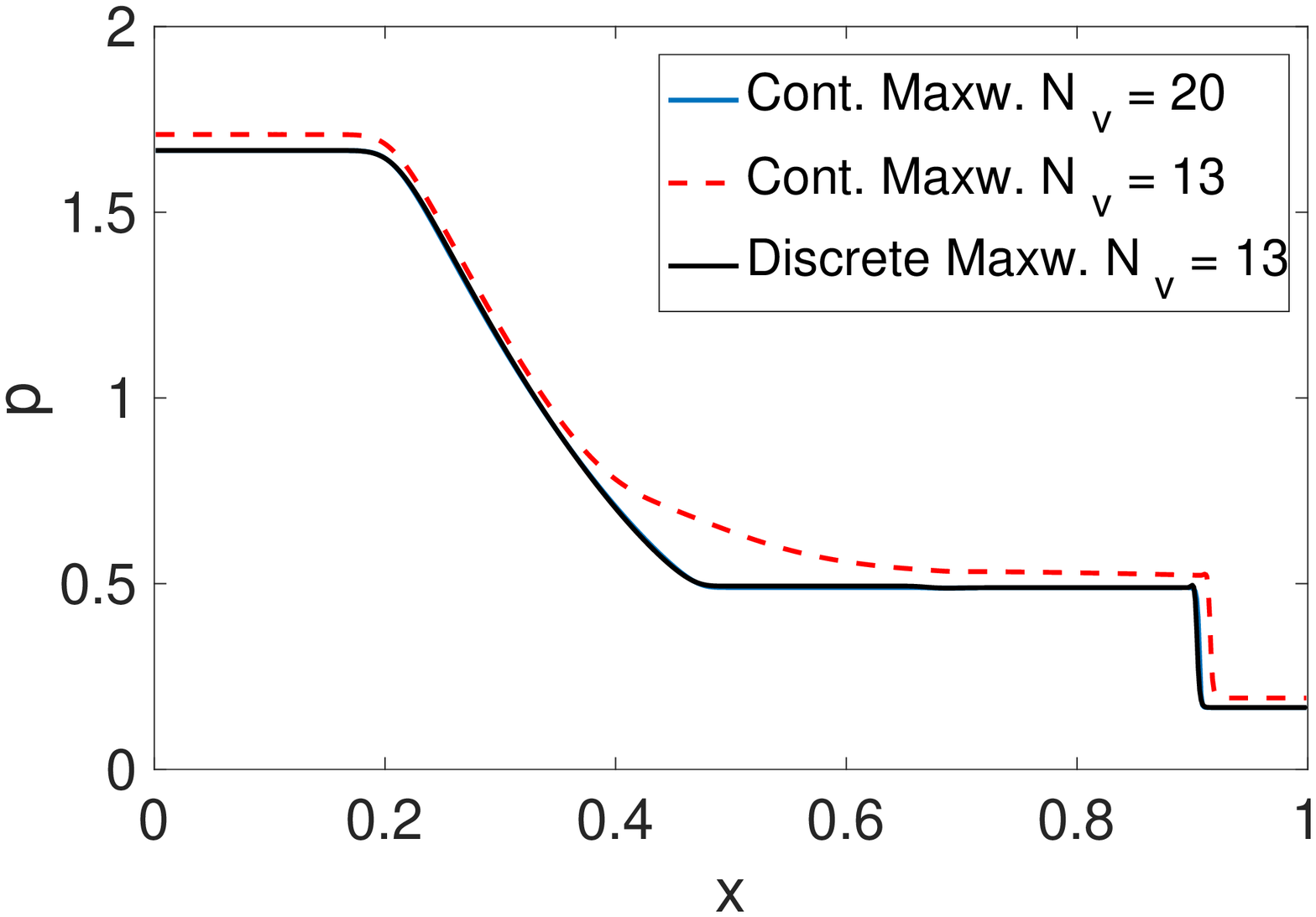}
\caption{Comparison of density, velocity and pressure from the  BGK model considering continuous Maxwellian  with $N_v=20$ and $N_v = 17$ together with the solutions considering discrete Maxwellian  with  $N_v=13$, $Nx = 800$, CFL = 1, $\lambda_l = 10^{-7}$ m for $x < 0.5$ , $\lambda_r =  8\times 10^{-7}$ m for $x > 0.5$ and $t_{\rm final } = 0.17$.}
\label{discrete_maxw}
\end{figure}

 In the case of three dimensional velocity space, the main drawback of the method is the memory problem as well as long computational time. Therefore, it is  important to reduce the computational time as well as memory allocation. Thus we have considered the standard Maxwellian as well as discrete Maxwellian.  
 If the number of velocity grid points is equal to $20$ or above, we do not see much difference of the solutions obtained from the standard as well as the discrete Maxwellian. However, the discrete Maxwellian requires 
 more computational efforts  than the standard Maxwellain because of the iterative method for solving nonlinear system of equations. If we choose the proper initial guess, the Newton iteration converges very fast. 
 We observed that if we consider discrete Maxwellian with $13$ velocity grids, the solutions are the same as the ones obtained from the standard Maxwellian with $20$ velocity grids. But with $13$ velocity grids with standard Maxwellian, the solutions deviate from the reference solutions. 
 In Figure \ref{discrete_maxw} we have plotted the density, velocity and pressure obtained from the BGK model with discrete Maxwellian considering 
 $N_v = 13$   and with the standard Maxwellian considering $N_v = 13$ together with  the reference solutions. We see that 
 the solutions obtained from the standard Maxwellian with $N_v =13$ deviates from the reference solutions, while the solutions with the discrete Maxwellain with $N_v = 13$ are closer with the reference solutions (standard Maxwellian with $N_v = 20$).   
 The computation is performed in $dual Intel Xeon Gold 6132 ("Skylake") @ 2.6 GHz $ with intel fortran compiler. 
 The total CPU to time with the standard Maxwellian having $20$ velocity grids is $443$ seconds  while with the discrete Maxwellian with $13$ velocity grids is $166$ seconds.   
 In higher dimensional cases the reduction of velocity grid points $N_v$ is very important from the computation as well as memory allocation point of view. 
 
\section{Conclusion and Outlook}
\label{sec:conclusion}

In this paper we have focussed on a meshfree method in the Semi Lagrangian scheme for the BGK model for rarefied gas flows. The meshfree method is applied for the reconstruction steps as well as for the 
implementation of boundary conditions. The diffuse boundary conditions on the solid wall is applied. The meshfree approximation is based on the moving least squares (MLS) method. 
We have used a linear approximation. 
The advantage of the 
meshfree approximation is that we do not require regular distribution of grid points in the velocity space, which will be very important if the boundary moves in time or interface between 
gases and other medium changes in time. Another advantage of the meshfree method is that we do not need to add ghost points. We have presented also the linear piecewise interpolations in order 
to compare the results from the linear MLS interpolations. We observed that linear MLS gives better results than the linear spline interpolation. We have considered a problem in one dimensional physical space and 
three dimensional velocity space. 
Sod's shock tube problem is solved for several ranges of mean free paths. 
No difference is found in the results obtained from regular and irregular grids in physical space. For larger mean free paths the solutions of BGK model are compared with the solutions obtained from the DSMC method for the Boltzmann equation. The solutions obtained from linear MLS are closer to DSMC results than the linear spline interpolation. For a very small mean free path we have compared the solutions of BGK model with the 
exact solutions of the compressible Euler equations and the solutions have very good agreement. In this case also the linear MLS approximation gives better results than the  linear spline interpolation. Moreover, we have shown that there is no restriction of the CFL number, where the CFL number can be larger than 1. 
Furthermore, we have studied the difference between the solutions obtained from the continuous and the discrete Maxwellian distribution. We found that the use of discrete Maxwellian allows us to reduce the number of grid points 
in velocity space without losing the accuracy in the solutions. The reduction of number of grid points is very important in higher dimensional physical spaces from the point of view of 
computational time and memory allocations. 

Future works will be the extension of the method for higher order reconstruction using WENO. Moreover, we are planning the extension of the method 
in higher dimensional physical spaces as well as the interaction of moving nano rigid particles immersed in a rarefied gas. 
  
\subsection*{Acknowledgment}  This  work is supported 
by the DFG (German research foundation) under Grant No. KL 1105/30-1 and by the ITN-ETN Marie-Curie Horizon 2020 program ModCompShock, Modeling and computation of shocks and interfaces, Project ID: 642768.

%
\end{document}